\documentclass[10pt]{article}%
\usepackage{amsmath}
\usepackage{graphicx}
\usepackage{amsfonts}
\usepackage{amssymb}%

\begin{document}
\noindent\textbf{\textsf{{\Large Studies of perturbed three vortex dynamics}}}

\bigskip

\quad\textsf{Denis Blackmore}

\quad\textit{Department of Mathematical Sciences, New Jersey Institute of
Technology, Newark, NJ 07102-1982}

\medskip

\quad\textsf{Lu Ting}

\quad\textit{Courant Institute of Mathematical Sciences, New York University,
New York, NY 10012-1186}

\medskip

\quad\textsf{Omar Knio}

\quad\textit{Department of Mechanical Engineering, The Johns Hopkins
University, Baltimore, MD 21218}

\bigskip

\begin{itemize}
\item[.] It is well known that the dynamics of three point vortices moving in
an ideal fluid in the plane can be expressed in Hamiltonian form,
where the resulting equations of motion are completely integrable in
the sense of Liouville and Arnold. The focus of this investigation
is on the persistence of regular behavior (especially periodic
motion) associated to completely integrable systems for certain
(admissible) kinds of Hamiltonian perturbations of the three vortex
system in a plane. After a brief survey of the dynamics of the
integrable planar three vortex system, it is shown that the
admissible class of perturbed systems is broad enough to include
three vortices in a half-plane, three coaxial slender vortex rings
in three-space, and `restricted' four vortex dynamics in a plane.
Included are two basic categories of results for admissible
perturbations: (i) general theorems for the persistence of invariant
tori and periodic orbits using Kolmogorov-Arnold-Moser and
Poincar\'{e}-Birkhoff type arguments; and (ii) more specific and
quantitative conclusions of a classical perturbation theory nature
guaranteeing the existence of periodic orbits of the perturbed
system close to cycles of the unperturbed system, which occur in
abundance near centers. In addition, several numerical simulations
are provided to illustrate the validity of the theorems as well as
indicating their limitations as manifested by transitions to chaotic
dynamics.
\end{itemize}

\ \ \ \textbf{Physics and Astronomy Classification Scheme:} 47.10.Df,
47.10.Fg, 47.15.ki, 47.20.Ky, 47.32.C-,

\ \ \ 47.32.cb, 47.52.+j

\vspace{0.3in}

\noindent\textsf{\textbf{I. INTRODUCTION}}

\medskip

Although virtually all research on vortex dominated fluid flows has
its roots in the seminal work of Helmholtz \cite{helm} (\emph{cf}.
\cite{lam}), specific advances in the dynamics of point vortices
moving in an \emph{ideal} (= inviscid, incompressible) fluid in a
planar region , which we shall refer to as $n$ \emph{vortex
dynamics} or the $n$ \emph{vortex problem}, can be traced back to
the pioneering work of Kirchhoff \cite{kirc} and Gr\"{o}bli
\cite{grob}. Kirchhoff was the first to describe the Hamiltonian
structure of the $n$ vortex problem, which he used to derive some
fundamental integrals of the motion, while Gr\"{o}bli conducted a
detailed analysis of the three vortex problem that included what was
essentially a proof of the integrability by quadratures of the
problem without employing Hamiltonian formalism, although he did not
give a complete description of the dynamics. About twenty years
after Gr\"{o}bli's remarkable work, Goryachev \cite{gory} took up
the problem again, and was able to obtain new insights concerning
the dynamics of point vortices. Shortly thereafter, Poincar\'{e}
\cite{poin} put his own imprimatur on vortex dynamics, just as he
did in so many other fields of research.

Some seventy years after the pioneering work of Kirchhoff and
Gr\"{o}bli, Synge \cite{syn} was able - using trilinear coordinates
- to fill many of the gaps in the dynamical picture of the three
vortex problem left by earlier studies. Among Synge's most important
contributions were the derivation of integrals and the
characterization of the critical points in terms of trilinear
coordinates, and the identification of a single parameter -
involving the sum of the product of the vortex strengths - that
distinguishes three distinct types of qualitative dynamics (the
elliptic, hyperbolic, and parabolic types) depending on whether this
parameter is positive, negative or zero. Twenty years later, Aref
\cite{aref, aref2} used the Hamiltonian based investigations of
Novikov \cite{nov} to rediscover the advantages of studying the
three vortex problem in the context of trilinear coordinates. In the
process, Aref identified additional properties of three vortex
dynamics, and initiated research on chaos in point vortex dynamics
and its relation to turbulent flows.

Subsequently, Tavantzis \& Ting \cite{tati}, taking Synge's approach
as their point of departure, derived a new constant of motion in
trilinear coordinates, which they used together with some classical
perturbation techniques to nearly complete the characterization of
three point vortex dynamics. Using linear analysis, they determined
the stability of all isolated stationary points in the trilinear
plane, and showed that expanding and contracting configurations are,
respectively, stable and unstable, which implies that contracting
similarity solutions and the eventual collision of three vortices is
unstable. More recently, Ting \emph{et al.}\cite{tbk} employed
techniques from nonlinear stability analysis to supply the few
missing details in the dynamics; in particular, they showed that
orbits starting just off the contracting configuration branch of the
singular curve of critical points in the parabolic case are
ultimately attracted to the expanding branch of this curve.

The Hamiltonian approach introduced by Kirchhoff \cite{kirc}, further
developed by Lin \cite{lin}, and perfected for point vortex problems by
Novikov \cite{nov} has proven to be very useful in vortex dynamics research.
Although not essential for solving the three vortex problem - as demonstrated
in \cite{grob}, \cite{syn}, \cite{tati}, and \cite{tbk} - presumably one could
use the integrals in involution to reduce it to a solvable
one-degree-of-freedom Hamiltonian system along the lines indicated by Borisov
and his collaborators in such papers as \cite{borleb}, \cite{bormamkil}, and
\cite{bormamkil2}. On the other hand, the symplectic structure underlying
Hamiltonian dynamics has proven to be extraordinarily effective in resolving a
wide variety of vortex problems such as the formulation of point vortex
dynamics on the sphere by Bogomolov \cite{bog}, the proof of complete
integrability of the three vortex problem on the sphere by Kidambi \& Newton
\cite{kidnew1} (see also \cite{borleb}, \cite{bormamkil2}, \cite{kidnew2} and
\cite{new}), verification of non-integrability of the general $n$ vortex
problem and $(n-1)$ coaxial vortex ring problem for $n>3$ by Bagrets \&
Bagrets \cite{bb} (see also Ziglin \cite{ziglin}), and several other important
results such as in \cite{aref2}, \cite{arefstrem}, \cite{eck}, \cite{kuzzas},
\cite{lim2}, \cite{sakajo}, and \cite{tok}.

It is in the Hamiltonian perturbation of integrable point vortex
dynamics where the Hamiltonian approach has proven to be
particularly useful. This is largely due to the availability of two
of the most important results in finite-dimensional symplectic
dynamics: Kolmogorov-Arnold-Moser (KAM) theory (such as in
\cite{arn1}, \cite{arn3}, \cite{kathas}, and \cite{mos1}) and
Poincar\'{e}-Birkhoff (PB) theory and its extensions and variants
(see \emph{e.g}. \cite{birk}, \cite{blwa}, \cite{bcw}, \cite{car},
\cite{coze}, \cite{ding}, \cite{flo}, \cite{frk}, \cite{gole},
\cite{jos}, \cite{kathas}, \cite{mos2}, \cite{henri}, and
\cite{rob}). Examples of such applications manifold: Khanin
\cite{khan} used a KAM theory inspired method to show that there
exist subsets of initial configurations for the (non-integrable)
four vortex problem in the plane leading to regular (integrable
like) motion; namely, quasiperiodic orbits on (deformed) invariant
(KAM) tori (\emph{cf}. Celletti \& Falcolini \cite{celfal}). A KAM
theory based argument combined with a deft application of Jacobi
canonical transformations enabled Lim \cite{lim1} to prove the
existence of quasiperiodic flow regimes in lattice vortex systems.
Blackmore \& Knio \cite{blkn}, using the fact that slender coaxial
vortex ring dynamics can be viewed as a perturbation of point vortex
motion, employed an innovative KAM theory type result to prove the
persistence of KAM tori and periodic orbits for an ample set of
initial positions of three rings sufficiently close to one another
when the rings having vortex strengths of the same sign. Later,
Blackmore \emph{et al}. \cite{bcw}, employing an analogous KAM
theory approach in concert with a novel extension of PB theory,
generalized this to any finite number of coaxial vortex rings with
strengths of the same sign. Similar results were obtained by
Blackmore \& Champanerkar \cite{blch} using the same type of
approach for any finite number of point vortices in a plane or
half-plane. In addition, Ting \& Blackmore \cite{tb} and Blackmore
\emph{et al}. \cite{btk} employed analogous ideas, together with
trilinear coordinates, to prove the persistence of regular flow
regimes for three coaxial vortex rings and three vortices in a
half-plane, respectively, when the vortex strengths are of differing
signs.

In this paper we shall, after a brief description of three vortex dynamics,
assemble and extend some of our recent results on the persistence of
quasiperiodic flows on KAM tori and periodic orbits for certain types of
non-integrable Hamiltonian perturbations of three vortex dynamics. Although we
shall include some interesting - fundamentally qualitative - results on the
existence of quasiperiodic and periodic orbits for perturbations of three
vortex dynamics, we intend to focus our attention on more quantitative
behavior concerning the persistence of periodic orbits that are, in an
appropriate sense, near periodic orbits of the unperturbed completely
integrable system. We shall obtain our results using a judicious combination
of Hamiltonian techniques and more classical perturbation methods that are
closely linked with the representation of the unperturbed dynamics in
trilinear coordinates.

The equations of motion for three vortex dynamics in both
Hamiltonian form and planar trilinear coordinates are presented in
Section \textsf{II}, along with a brief review of the trilinear
phase portrait that emphasizes those dynamical features that
exploited extensively in the sequel. Next, in Section \textsf{III},
we describe specific types of perturbations that are subsumed by the
Hamiltonian perturbations for which our main conclusions are proven
in subsequent sections. The specific perturbations considered are
three vortices in a half-plane, a `restricted' four vortex problem
in the plane, and three coaxial vortex rings. We prove qualitative
theorems in Section \textsf{IV} on the persistence of quasiperiodic
flows on KAM tori interspersed with periodic orbits for Hamiltonian
perturbations of three vortex dynamics subject to hypotheses
satisfied by the examples in the preceding section. A more classical
perturbation approach is used in Section \textsf{V} to show that the
perturbations of type introduced in Section \textsf{IV} have the
property of having periodic orbits that are close - in an
appropriate coordinate system - to certain cycles of the three
vortex system.

In the penultimate part of the paper, Section \textsf{VI}, we
present a variety of numerical simulations that illustrate the
persistence of regularity associated with the three vortex problem
under perturbations satisfying the properties in our main results.
The properties in question provide sufficient conditions on the
initial configurations of point vortices or slender coaxial vortex
rings for the existence of invariant tori and periodic orbits, and
for periodic orbits of the perturbed system that are close to those
of the unperturbed system. A number of simulations are also included
to show how regularity breaks down - signalled by the emergence of
chaos - as the limitations imposed by the hypotheses in our theorems
are exceeded. Finally, in Section \textsf{VII}, we summarize and
distill the most important conclusions in the paper, and indicate
some promising directions for related future research.

\bigskip

\noindent\textsf{\textbf{II. UNPERTURBED DYNAMICS}}

\medskip

In this section we describe the equations of motion of point vortices in an
ideal fluid in the plane, and give a rather complete characterization of the
dynamics for the integrable three vortex system. Knowing the dynamics in the
three vortex case shall prove quite useful in our subsequent analysis and
description of periodic orbits for perturbed systems that are perturbations of
cycles of the unperturbed system.

\medskip

\noindent\textsf{\textbf{A. Governing equations}}

\smallskip

We begin with the general problem of $n$ point vortices of respective nonzero
strengths $\Gamma_{1},...,\Gamma_{n}$ moving in an ideal fluid in the
(complex) plane $\mathbb{C}$ $(=\mathbb{R}^{2})$ and located at the positions
$z_{k}=x_{k}+iy_{k}$, $1\leq k\leq n$, respectively. The equations of motion
in complex form are%
\begin{equation}
\dot{\bar{z}}_{k}=i%
\sum\limits_{j=1,j\neq k}^{n}
\kappa_{j}\left(  z_{j}-z_{k}\right)  ^{-1},\qquad(1\leq k\leq n) \label{eq1}%
\end{equation}
where the overbar indicates the complex conjugate, the overdot
denotes differentiation with respect to $t$, and
$\kappa_{j}:=\Gamma_{j}/2\pi$, $1\leq
j\leq n$, which can be recast as the complex Hamiltonian equation%
\begin{equation}
\mathbf{\kappa}\ast\dot{\bar{\mathbf{z}}}=2i\partial_{\mathbf{z}}H_{0},
\label{eq2}%
\end{equation}
where $\mathbf{\kappa}:=\left(  \kappa_{1},...,\kappa_{n}\right)  $,
$\mathbf{z}:=\left(  z_{1},...,z_{n}\right)  $, $\ast$ denotes the usual
(Cartesian product ring) product in $\mathbb{C}^{n}$, $\partial_{\mathbf{z}%
}H_{0}:=\left(  \partial_{z_{1}}H_{0},...,\partial_{z_{n}}H_{0}\right)  $, and
the Hamiltonian function is given as
\begin{equation}
H_{0}:=- {\sum\limits_{1\leq j<k\leq n}}
\kappa_{j}\kappa_{k}\log\left\vert z_{j}-z_{k}\right\vert . \label{eq3}%
\end{equation}
We note here that to guarantee smoothness of the system (1), the phase space
needs to be defined as%
\[
\mathbb{C}_{\#}^{n}:=\left\{  \mathbf{z}\in\mathbb{C}^{n}:z_{j}\neq
z_{k}\;\forall j\neq k\right\}  ,
\]
and we introduce the following notation that will prove useful in the sequel:%
\[
\mathfrak{K}_{n}^{(1)}:=%
\sum\limits_{j=1}^{n}
\kappa_{j},\quad\mathfrak{K}_{n}^{(2)}:=%
{\sum\limits_{1\leq j<k\leq n}} \kappa_{j}\kappa_{k}.
\]
The governing equations can be converted into the more familiar real
Hamiltonian form in $\mathbb{R}^{2n}$ (or more properly $\mathbb{R}_{\#}%
^{2n})$
\begin{equation}
\dot{x}_{k}=\kappa_{k}^{-1}\partial_{y_{k}}H_{0}=\left\{  H_{0},x_{k}\right\}
,\quad\dot{y}_{k}=-\kappa_{k}^{-1}\partial_{x_{k}}H_{0}=\left\{  H_{0}%
,y_{k}\right\}  ,\qquad(1\leq k\leq n) \label{eq4}%
\end{equation}
where $\mathbb{R}_{\#}^{2n}$ is $\mathbb{C}_{\#}^{n}$ in the usual
representation in real coordinates, and the (nonstandard) Poisson bracket is
defined as%
\[
\left\{  f,g\right\}  :=%
\sum\limits_{k=1}^{n} \kappa_{k}^{-1}\left(  \frac{\partial
f}{\partial y_{k}}\frac{\partial g}{\partial x_{k}}-\frac{\partial
f}{\partial x_{k}}\frac{\partial g}{\partial y_{k}}\right)  .
\]

It is easy to see that
\begin{equation}
H_{0},\;J:=%
\sum\limits_{k=1}^{n}
\kappa_{k}z_{k},\;K:=%
\sum\limits_{k=1}^{n}
\kappa_{k}\left\vert z_{k}\right\vert ^{2}, \label{eq5}%
\end{equation}
representing the total energy, the linear momentum (impulse) and the angular
momentum, respectively, are integrals of the system (1) (= (2) = (4)). From
these we can construct the following three functionally independent, real
constants of motion that are in involution:
\begin{equation}
H_{0},\;K,\;L:=\left(  \mathfrak{R}J\right)  ^{2}+\left(  \mathfrak{I}%
J\right)  ^{2}, \label{eq6}%
\end{equation}
where $\mathfrak{R}$ and $\mathfrak{J}$denote the real and imaginary part,
respectively, of a complex number, and it is easy to verify the involutivity
of these integrals by computing that%
\begin{equation}
\left\{  H_{0},K\right\}  =\left\{  H_{0},L\right\}  =\left\{  K,L\right\}
=0. \label{eq7}%
\end{equation}
Accordingly the dynamics of three point vortices in a plane is completely
integrable in the sense of Liouville and Arnold, which we shall denote as
\emph{LA-integrable}. Whence it follows from Liouville-Arnold theory and the
Poincar\'{e}-Birkhoff \ fixed point theorem and its extensions (see
\emph{e.g}. \cite{arn3}, \cite{birk}, \cite{blwa}, \cite{blch}, \cite{car},
\cite{coze}, \cite{ding}, \cite{flo}, \cite{frk}, \cite{gole}, \cite{jos},
\cite{kathas}, \cite{new}, and \cite{henri}) that three vortex dynamics is
quasiperiodic, confined to invariant 3-tori and exhibits periodic orbits of
arbitrarily large periods for all combinations of nonvanishing vortex strengths.

\medskip

\noindent\textsf{\textbf{B. Dynamics in trilinear coordinates}}

\smallskip

Three vortex dynamics, which we summarize here, can be described in a very
efficient manner using trilinear coordinates, which have their roots in
algebraic geometry.. Our approach follows that of Synge \cite{syn}, Tavantzis
\& Ting \cite{tati}, and Ting \emph{et al. }\cite{tbk}, and we refer the
reader to these sources and Aref \cite{aref, aref2} for further details. We
begin with a dynamic formulation introduced by Gr\"{o}bli \cite{grob} in terms
of the sides of the (possibly degenerate) triangular configuration denoted as
$R_{1}:=\left\vert z_{2}-z_{3}\right\vert $, $R_{2}:=\left\vert z_{1}%
-z_{3}\right\vert $, and $R_{3}:=\left\vert z_{1}-z_{2}\right\vert $, which
can be expressed as%
\begin{align}
R_{1}\dot{R}_{1}  &  =2A\kappa_{1}\left(  R_{2}^{-1}-R_{3}^{-1}\right)
,\nonumber\\
R_{2}\dot{R}_{2}  &  =2A\kappa_{2}\left(  R_{3}^{-1}-R_{1}^{-1}\right)
,\label{eq8}\\
R_{3}\dot{R}_{3}  &  =2A\kappa_{3}\left(  R_{1}^{-1}-R_{2}^{-1}\right)
,\nonumber
\end{align}
where $A$ denotes the oriented area of the configuration defined to be
positive (negative) for a counterclockwise (clockwise) arrangement of the
vertices ordered as $z_{1},z_{2},z_{3}$. Naturally the points $\mathbf{R}%
=(R_{1},R_{2},R_{3})$ in Euclidean three-space that comprise the domain of (8)
must be confined to the first octant (and also avoid or compensate for
singularities on the bounding coordinate planes). It is interesting to note
that the stationary (fixed) points of (8) correspond to equilateral
configurations, or collinear configurations satisfying the additional
condition $\dot{A}=0$, where the extra condition in the collinear case is
necessitated by the singular nature of the vector field at such points.
Observe that the systems (1), (4) and (8) are invariant under the
transformation $t\rightarrow-t,\kappa_{j}\rightarrow-\kappa_{j}\;(1\leq
j\leq3)$, hence we may assume without loss of generality that%
\begin{equation}
\kappa_{1}\geq\kappa_{2}>0\text{ and }\kappa_{2}\geq\kappa_{3}, \label{e9}%
\end{equation}
which we shall do throughout the sequel.

The triangle inequality imposes additional restrictions on the
domain of (8). We shall define this restricted domain in accordance
with the approach of Synge \cite{syn} (and Tavantzis \& Ting
\cite{tati}), wherein one can find some excellent figures
representing the concepts described below, by first introducing the
equilateral triangle $\mathcal{T}:=$ $P_{1}P_{2}P_{3}$ of height one
determined by intersecting the first octant with the plane
$\mathcal{P}$ defined as the locus of
$R_{1}+R_{2}+R_{3}=\sqrt{2/3}$, where $P_{1}:=\left(
\sqrt{2/3},0,0\right)  $, $P_{2}:=\left(  0,\sqrt {2/3},0\right)  $,
and $P_{3}:=\left(  0,0,\sqrt{2/3}\right)  $. Observe that the first
octant can be viewed as the fiber bundle over $\mathcal{T}$ with the
rays emanating from the origin (but not including the origin) as
fibers. It is easy to verify from the triangle inequality that the
proper base space of the bundle defining the admissible points
$\mathcal{D}$ for (8) is the subtriangle $T=Q_{1}Q_{2}Q_{3}$, where
$Q_{1}$ is the midpoint of the edge $P_{2}P_{3}$, $Q_{2}$ is the
midpoint of the edge $P_{1}P_{3}$, and $Q_{3}$ is the midpoint of
the edge $P_{1}P_{2}$, so that $Q_{1}:=\left(  0,\sqrt{1/6},\sqrt
{1/6}\right)  $, $Q_{2}:=\left(  \sqrt{1/6},0,\sqrt{1/6}\right) $,
and $Q_{3}:=\left(  \sqrt{1/6},\sqrt{1/6},0\right)  $. We note that
the bundle projection $\rho:\mathcal{D}\rightarrow T$ is just the
radial projection on $T$.

Points in the triangle $T$ can be assigned \emph{trilinear coordinates}
$\mathbf{x}=\left(  x_{1},x_{2},x_{3}\right)  $ defined as%
\begin{equation}
x_{j}:=R_{j}\left(  R_{1}+R_{2}+R_{3}\right)  ^{-1},\quad(1\leq j\leq3)
\label{e10}%
\end{equation}
which we note are not independent inasmuch as they satisfy $x_{1}+x_{2}%
+x_{3}=1$, and represent distance from the edges of $\mathcal{T}$ ; in
particular, $x_{1}$, $x_{2}$, and $x_{3}$ are, respectively, the distances of
$\mathbf{x}$ from the edges $P_{2}P_{3}$, $P_{1}P_{3}$, and $P_{1}P_{2}$. The
points in the planar triangle $T$ can also be conveniently represented in
terms of the following coordinates:%
\begin{equation}
\alpha:=\left(  1/\sqrt{3}\right)  \left(  x_{2}-x_{1}\right)  ,\quad
\beta:=x_{3}, \label{e11}%
\end{equation}
with `inverse' transformation
\begin{equation}
x_{1}=\left(  1/2\right)  \left(  1-\beta-\alpha\sqrt{3}\right)
,\;x_{2}=\left(  1/2\right)  \left(  1-\beta+\alpha\sqrt{3}\right)
,\;x_{3}=\beta. \label{e12}%
\end{equation}
The system (8) projected on $T$ assumes the following form in trilinear
coordinates%
\begin{align}
\dot{x}_{1}  &  =\Lambda\left\{  \kappa_{1}x_{1}\left(  x_{3}^{2}-x_{2}%
^{2}\right)  -x_{1}\left[  \kappa_{1}x_{1}\left(  x_{3}^{2}-x_{2}^{2}\right)
+\kappa_{2}x_{2}\left(  x_{1}^{2}-x_{3}^{2}\right)  +\kappa_{3}x_{3}\left(
x_{2}^{2}-x_{1}^{2}\right)  \right]  \right\}  ,\nonumber\\
\dot{x}_{2}  &  =\Lambda\left\{  \kappa_{2}x_{2}\left(  x_{1}^{2}-x_{3}%
^{2}\right)  -x_{2}\left[  \kappa_{1}x_{1}\left(  x_{3}^{2}-x_{2}^{2}\right)
+\kappa_{2}x_{2}\left(  x_{1}^{2}-x_{3}^{2}\right)  +\kappa_{3}x_{3}\left(
x_{2}^{2}-x_{1}^{2}\right)  \right]  \right\}  ,\label{e13}\\
\dot{x}_{3}  &  =\Lambda\left\{  \kappa_{3}x_{3}\left(  x_{2}^{2}-x_{1}%
^{2}\right)  -x_{3}\left[  \kappa_{1}x_{1}\left(  x_{3}^{2}-x_{2}^{2}\right)
+\kappa_{2}x_{2}\left(  x_{1}^{2}-x_{3}^{2}\right)  +\kappa_{3}x_{3}\left(
x_{2}^{2}-x_{1}^{2}\right)  \right]  \right\}  ,\nonumber
\end{align}
where
\begin{equation}
\Lambda:=2A\left[  \left(  R_{1}R_{2}R_{3}\right)  ^{-1}\left(  R_{1}%
+R_{2}+R_{3}\right)  \right]  ^{2}. \label{e14}%
\end{equation}
The dimensionality of (12) can be reduced by one by recasting the system in
$\alpha,\beta$ - coordinates as%
\begin{align}
\dot{\alpha}  &  =\Lambda\left\{  \left(  1/8\sqrt{3}\right)  \left[
\kappa_{2}(1-3\beta-\alpha\sqrt{3})\left(  1-(\beta-\alpha\sqrt{3}%
)^{2}\right)  +\right.  \right. \nonumber\\
&  \left.  \left.  \quad\kappa_{1}(1-3\beta+\alpha\sqrt{3})\left(
1-(\beta+\alpha\sqrt{3})^{2}\right)  \right]  -\alpha\Xi\right\}
,\label{e15}\\
\dot{\beta}  &  =\Lambda\left\{  \kappa_{3}\sqrt{3}\alpha\beta\left(
1-\beta\right)  -\left(  \beta/8\right)  \Xi\right\}  ,\nonumber
\end{align}
where%
\begin{align}
\Xi &  :=\left\{  -\kappa_{1}\left(  1-3\beta+\alpha\sqrt{3}\right)  \left[
1-\left(  \beta+\alpha\sqrt{3}\right)  ^{2}\right]  +\right. \nonumber\\
&  \left.  \quad\kappa_{2}\left(  1-3\beta-\alpha\sqrt{3}\right)  \left[
1-\left(  \beta-\alpha\sqrt{3}\right)  ^{2}\right]  +\kappa_{3}8\sqrt{3}%
\alpha\beta\left(  1-\beta\right)  \right\}  , \label{e16}%
\end{align}
and the factor $\Lambda$ can be absorbed into a rescaling of time in order to
simplify the description of the phase portrait on $T$.

The system (13) or (15) actually includes two possible orientations
of the triangular configurations of three vortices, so along with
the phase portrait on $T$, we must include another copy of $T$,
which we denote as $T^{\ast}$, having the opposite orientation. More
precisely, the integral curves on $T^{\ast}$ are precisely those of
$T$ except they have the opposite orientation with respect to $t$.
By gluing $T$ and $T^{\ast}$ (and their vector fields) smoothly
together along their corresponding sides to obtain the double
$M:=D(T)$ of $T$, which is a 2-sphere, we obtain a vector field $X$
on $M$ by joining the oriented field corresponding to (13) or (15).
Note that at this point, we can reformulate the bundle structure
described above for $\mathcal{D}$ in the `doubled' form
$\rho:\mathcal{D}\rightarrow M=D(T)$. As for fixed points of the
flow, we always have, in terms of trilinear coordinates,
$E:=(1/3,1/3,1/3)$ on $T$ and its copy $E^{\ast}$ on the $T^{\ast}$
half of $M$ corresponding to an equilateral configuration, and the
common vertices $Q_{1}=Q_{1}^{\ast}$, $Q_{2}=Q_{2}^{\ast}$ and $Q_{3}%
=Q_{3}^{\ast}$ corresponding to points where a pair of vortices coincide, as
well as additional points that we shall describe in what follows.

A local stability analysis at $E$ (or $E^{\ast}$) shows that there are three
distinct types of dynamical behavior corresponding to the value of
$\mathfrak{K}^{(2)}:=\mathfrak{K}_{3}^{(2)}$: These are the \emph{elliptic}
case when $\mathfrak{K}^{(2)}>0$, the \emph{hyperbolic} case when
$\mathfrak{K}^{(2)}<0$, and the exceptional \emph{parabolic} case when
$\mathfrak{K}^{(2)}=0$. When the system is elliptic, $E$ (and $E^{\ast}$) is a
center; $E$ (and $E^{\ast}$) is a saddle point in the hyperbolic case; and
when the system is parabolic, $E$ (and $E^{\ast}$) is no longer isolated - it
lies on the curve of fixed points $C$ on $T$ (with copy $C^{\ast}$ on
$T^{\ast}$) defined as%
\[
C:\kappa_{1}^{-1}x_{1}^{2}+\kappa_{2}^{-1}x_{2}^{2}+\kappa_{3}^{-1}x_{3}%
^{2}=0.
\]

One of the most remarkable and useful features of the trilinear coordinate
based projection of the system (8) on $M$ is what amounts to essentially a
unique path lifting property (cf. \cite{span}) for integral curves on $M$,
except along $C$ for parabolic systems. For future reference, we shall refer
to this property as the \emph{unique integral path lifting principle} for the
bundle $\rho:\mathcal{D}\rightarrow M=D(T)$. More precisely, the initial
points in $\mathcal{D}$ for integral curves of (8) establishes a bijective
correspondence via the projection $\rho$ with the integral curves of (12) or
(14) on $M$, except along $C$ in the parabolic case. This bijective
correspondence can be readily established using the integrals of motion of the
systems. As for these flow invariants, employing the constants of motion
$H_{0}$ and $K$ for (4), it is easy to verify that
\begin{equation}
\kappa_{1}^{-1}\log R_{1}+\kappa_{2}^{-1}\log R_{2}+\kappa_{3}^{-1}\log
R_{3}\text{ and }\kappa_{1}^{-1}R_{1}^{2}+\kappa_{2}^{-1}R_{2}^{2}+\kappa
_{3}^{-1}R_{3}^{2} \label{e17}%
\end{equation}
are first integrals for (8). These integrals were cleverly combined by
Tavantzis \& Ting \cite{tati} to obtain the following invariants for (13):%
\begin{equation}
I:=\left( {\sum\limits_{k=1}^{3}} \kappa_{k}^{-1}x_{k}^{2}\right)
\left( {\prod\limits_{k=1}^{3}}
x_{k}^{2/\kappa_{k}}\right)  ^{(\kappa_{1}\kappa_{2}\kappa_{3}/\mathfrak{K}%
^{(2)})}, \label{e18}%
\end{equation}
when $\mathfrak{K}^{(2)}\neq0$, and
\begin{equation}
I_{0}:=\left(  x_{1}/x_{3}\right)  ^{2\kappa_{2}}\left(  x_{2}/x_{3}\right)
^{2\kappa_{1}} \label{e19}%
\end{equation}
for the parabolic case. The above integrals can be combined (see \cite{btk}
and \cite{tbk}) to produce a an invariant valid for all cases, namely%
\begin{equation}
\bar{I}:=\left( {\sum\limits_{k=1}^{3}}
\kappa_{k}^{-1}x_{k}^{2}\right)
^{\mathfrak{K}^{(2)}/(2\kappa_{3})}\left(
x_{1}^{1/\kappa_{1}}x_{2}^{1/\kappa_{2}}x_{3}^{1/\kappa_{3}}\right)
^{\kappa_{1}\kappa_{2}}. \label{e20}%
\end{equation}
We note that we shall find the invariant $\bar{I}$ to be particularly useful
for our numerical simulations of perturbations of three vortex dynamics in the
sequel. In particular, the extent to which $\bar{I}$ differs from a constant
value for a perturbed system is a good indication of the extent to which the
dynamics of the systems diverges from the LA-integrable motion of the three
vortex problem. Any of the above constants of motion can be used to prove the
almost unique path lifting property (\emph{cf}. Synge \cite{syn}).

Before embarking on a condensed - but reasonably complete - trilinear
coordinate based description of the dynamics in the elliptic, hyperbolic and
parabolic cases, certain aspects of the compelling effectiveness of the
trilinear approach almost demand comment. We first observe that (8) - and so
also (13) - incorporates the Lie group of symmetries of the original system
(1) - which is the planar Euclidean group. As one would expect to be able to
use the symmetries in a standard symplectic fashion via the associated
momentum map to obtain a Hamiltonian restriction of the original system (4)
having just one degree of freedom for the three vortex problem, there is the
obvious suggestion of a strong link between the symplectic reduction and the
system (13). It appears that such a link, which shall have no need of in the
sequel, has yet to be definitively established.

\medskip

\noindent\textbf{Elliptic case} \textbf{(}$\mathbf{\mathfrak{K}}%
^{(2)}\mathbf{>0}$\textbf{) }

\smallskip

As indicated above, in the typical elliptic case, which is
illustrated for $T$ in Fig. 1, both the stationary point $E$ and its
(opposite orientation) copy $E^{\ast}$ are centers with index $+1$.
It can be easily shown that the fixed points $Q_{1}$, $Q_{2}$ and
$Q_{3}$ are also centers. In addition, there are three other fixed
points of (13) - one on the interior of each side of $T$ - denoted
as $Q_{4}$, $Q_{5}$ and $Q_{6}$ in Fig. 1, and each of these
additional stationary points can readily be shown to be a saddle
points of index $-1$. Notice that these observations are consistent
with the Poincar\'{e}-Hopf Index theorem (see \emph{e.g}.
\cite{kathas}) , inasmuch as
a simple calculation yields%
\[
5(+1)+3(-1)=2=\chi\left(  M\right)  =\chi\left(  D(T)\right)  ,
\]
where $\chi(M)$ is the Euler-Poincar\'{e} characteristic of the sphere $M$,
which is two. In this case, we see that the global phase portrait of (15) is
completely determined, and this leads to an exhaustive determination of three
vortex dynamics in the elliptic case owing to the unique integral path lifting
principle mentioned above.

\begin{figure}[ptbh]
\begin{center}
\includegraphics[height=2.6in.]{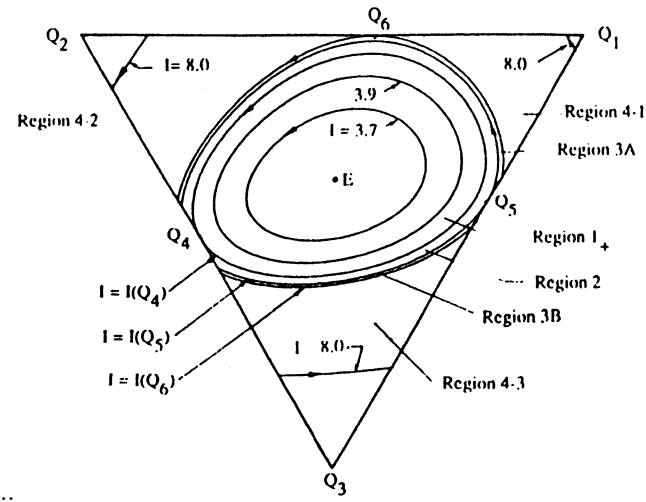}
\end{center}
\caption{Elliptic type with $\kappa_{1}=2$, $\kappa_{2}=1$, and $\kappa
_{3}=1/2.$}%
\end{figure}

\medskip

\noindent\textbf{Hyperbolic case} \textbf{(}$\mathbf{\mathfrak{K}}%
^{(2)}\mathbf{<0}$\textbf{)}

\smallskip

An example of the dynamics for the hyperbolic case is shown in Fig. 2. For all
hyperbolic cases, the points $E$ and $E^{\ast}$ associated to equilateral
configurations are stationary saddle points of index $-1$. When it comes to
the fixed points at the vertices of the triangle $T$ (glued to those of
$T^{\ast}$), the point $Q_{3}$ is always a center of index $+1$ for all
possible admissible vortex strengths. Any other fixed points, which of course
must include the other vertices of $T$ as well other points on the open edges
of the triangle can have a variety of natures depending on certain algebraic
criteria (see \cite{tati}). For example, in the case shown in Fig. 2, $Q_{1}$
and $Q_{2}$ are both centers, while the stationary points $Q_{4}$, $Q_{5}$ and
$Q_{6}$ are, respectively a center, center and saddle point. Thus we compute
the following index sum%
\[
\mathrm{ind}E+\mathrm{ind}E^{\ast}+\mathrm{ind}Q_{3}+\mathrm{ind}%
Q_{1}+\mathrm{ind}Q_{2}+\mathrm{ind}Q_{4}+\mathrm{ind}Q_{5}+\mathrm{ind}%
Q_{6}=2(-1)+5(+1)+(-1)=2,
\]
which is consistent with the Poincar\'{e}-Hopf Index theorem.

Depending on the algebraic criteria (involving the vortex strengths), there
are exceptional cases where $Q_{5}$ merges with $Q_{1}$ or $Q_{4}$ coincides
with $Q_{2}$, or possibly both, in which case one or both of $Q_{1}$ and
$Q_{2}$ can assume the form of a degenerate isolated stationary point of index
$+2$. In all cases, the corresponding index sum can be shown to be in
agreement with the Poincar\'{e}-Hopf formula. Exceptional cases
notwithstanding, once again we also have, for any parameter values consistent
with the hyperbolic case, a complete characterization of three vortex dynamics
owing to the unique integral path lifting principle.

\begin{figure}[ptbh]
\begin{center}
\includegraphics[height=2.7in.]{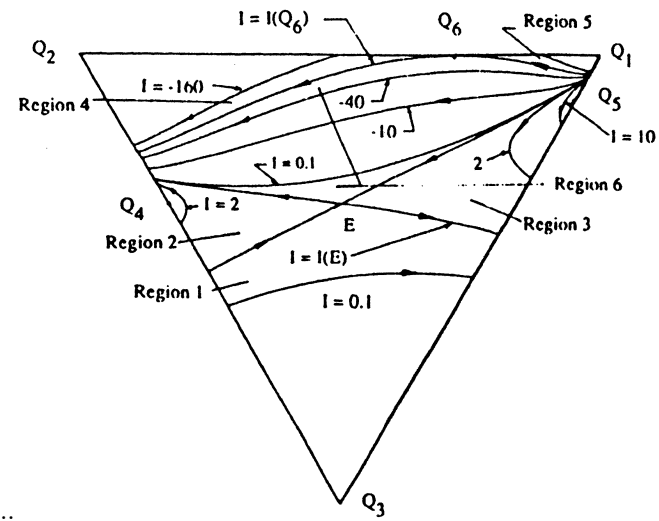}
\end{center}
\caption{Hyperbolic type with $\kappa_{1}=2$, $\kappa_{2}=1$, and $\kappa
_{3}=-\,4/5.$}%
\end{figure}

\begin{figure}[ptbh]
\begin{center}
\includegraphics[height=3.6in.]{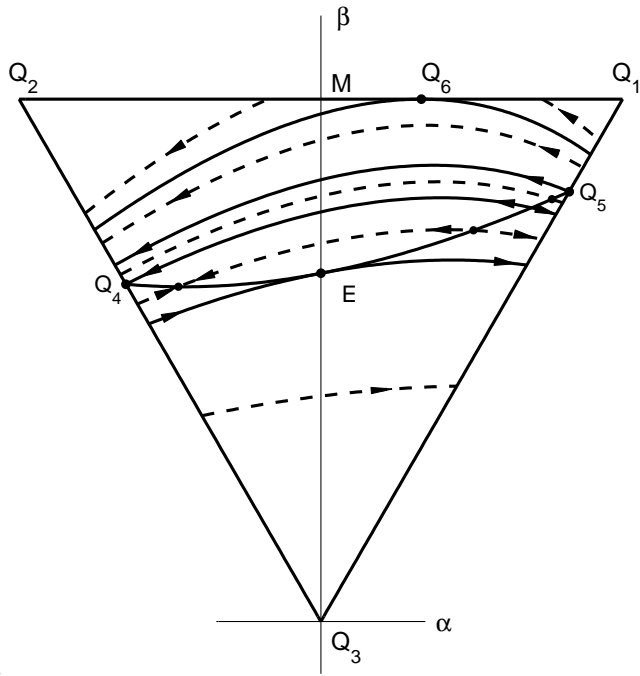}
\end{center}
\caption{Parabolic type with $\kappa_{1}=2$, $\kappa_{2}=1$, and
$\kappa
_{3}=-\,2/3.$}%
\end{figure}

\medskip

\noindent\textbf{Parabolic case} \textbf{(}$\mathbf{\mathfrak{K}}%
^{(2)}\mathbf{=0}$\textbf{)}

\smallskip

Three vortex dynamics represented in the trilinear coordinate phase plane for
an example of the parabolic case is illustrated in Fig. 3. Observe, as was
discussed above, in this case the points $E$ and $E^{\ast}$ lie on the
singular curve $C$ in $T$ and its copy $C^{\ast}$ on $T^{\ast}$ comprised of
stationary points of the system (13). Unfortunately, the unique integral path
lifting principle does not apply on the double $D\left(  C\right)  $ of $C$,
which is a circle on $M$ to which the points $E$ and $E^{\ast}$ both belong,
and with this is associated another problem with the singular curve exhibited
in Fig. 3; namely, the trajectories of (13) appear to cross this curve, which
is not what one expects for a reasonable dynamical system in which standard
uniqueness theorems obtain. But the principle does apply on the complement of
$D\left(  C\right)  $ on $M$, making it possible to obtain an almost complete
description of three vortex dynamics in the parabolic case. To fill in the
gaps, it is necessary only to conduct a more careful analysis of the dynamics
of the full system in a neighborhood of the singular circle $D\left(
C\right)  $ as in Ting \emph{et al}.\cite{tbk}, but we shall not pursue this
further here since it not needed for our analysis in the sequel.

Note that if we cut the sphere along the singular curve $D\left(  C\right)  $
in the example shown in Fig. 3, we obtain two disks $B_{1}$ and $B_{2}$ for
which we can apply the extension of the Poincar\'{e}-Hopf formula to manifolds
with boundaries. The upper disk $B_{1}$ contains the stationary points $Q_{1}%
$, $Q_{2}$ and $Q_{6}$ in its interior, and it is easy to verify that these
fixed points are, respectively, a center, center, and saddle point. Whence, we
find that the index sum in $B_{1}$ is%
\[
2(+1)+(-1)=1=\chi\left(  B_{1}\right)  .
\]
For the lower hemisphere $B_{2}$, we find that $Q_{3}$ is the only interior
fixed point, and it is a center of index $+1$. Accordingly we compute that, as
expected,
\[
1(+1)=1=\chi\left(  B_{2}\right)  .
\]

Depending on the parameter values of the system for a particular example of
the parabolic case, certain dynamic properties can vary from the example
depicted in Fig. 3. However, as shown in Tavantzis \& Ting \cite{tati} (and
can also be inferred from the Poincar\'{e}-Hopf formula in virtue of the fact
that $Q_{3}$ is the only stationary point below the singular curve $C$), the
fixed point $Q_{3}$ is always a center.

Before leaving our summary of the trilinear coordinate phase plane behavior of
the three vortex problem, we wish to emphasize the fact that regardless of the
type, the stationary point $Q_{3}$ is always a center, and there are often
additional centers.

\bigskip

\noindent\textsf{\textbf{III. TYPES OF PERTURBATIONS}}

\medskip

Here we shall provide examples of perturbations three point dynamics that
satisfy the hypotheses of the main results that we shall derive in succeeding
sections. These particular types of perturbations will also be used for
numerical illustrations of where the theorems apply, and where they begin to
break down as evidenced by the onset of nonregular, chaotic regimes. The types
of perturbations that we choose are three vortex dynamics in a half-plane,
restricted four vortex dynamics, and three coaxial vortex ring dynamics in space.

\medskip

\noindent\textsf{\textbf{A. Three vortex dynamics in a half-plane}}

\smallskip

The complex Hamiltonian governing equations for three point vortices of
nonzero strengths $\Gamma_{k}$ at respective points $z_{k}=x_{k}+iy_{k}$,
$1\leq k\leq3$, in motion in an ideal fluid in the half-plane $\mathbb{H}%
:=\{z=x+iy\in\mathbb{C}:y>0\}$ are%
\begin{equation}
\kappa_{k}\dot{\bar{z}}_{k}=i%
{\sum\limits_{j=1,j\neq k}^{3}}
\kappa_{j}\kappa_{k}\left(  z_{j}-z_{k}\right)  ^{-1}-i%
{\sum\limits_{j=1}^{3}}
\kappa_{j}\kappa_{k}\left(  \bar{z}_{j}-z_{k}\right)  ^{-1}=2i\partial_{z_{k}%
}H_{A},\qquad(1\leq k\leq3) \label{e21}%
\end{equation}
where the Hamiltonian function is%
\begin{equation}
H_{A}:=-%
{\sum\limits_{1\leq j<k\leq3}}
\kappa_{j}\kappa_{k}\log\left\vert z_{j}-z_{k}\right\vert +%
{\sum\limits_{1\leq j\leq k\leq3}}
\kappa_{j}\kappa_{k}\log\left\vert \bar{z}_{j}-z_{k}\right\vert . \label{e22}%
\end{equation}
With the understanding, expressed in Section \textsf{II}, that
smoothness requirements
actually demand that the domain be restricted to points in $\mathbb{H}%
_{\#}^{3}:=\{\mathbf{z}\in\mathbb{H}^{3}:z_{j}\neq z_{k}\;\forall j\neq k\}$,
we can express (21) in smooth, real Hamiltonian form%
\begin{equation}
\dot{x}_{k}=\kappa_{k}^{-1}\partial_{y_{k}}H_{A}=\left\{
H_{A},x_{k}\right\}
,\quad\dot{y}_{k}=-\kappa_{k}^{-1}\partial_{x_{k}}H_{A}=\left\{  H_{A}%
,y_{k}\right\}  ,\qquad(1\leq k\leq3) \label{e23}%
\end{equation}
where we use the same Poisson bracket as in (4). This system, as is easily
verified, has the following independent constants of motion in involution%
\begin{equation}
H_{A},\;\mathfrak{I}J:=\kappa_{1}y_{1}+\kappa_{2}y_{2}+\kappa_{3}y_{3},\;
\label{e24}%
\end{equation}
and this appears to be the maximal such set of invariants: Although we are
unaware of a proof showing that (21), or certain analogs for vortex dynamics
on a portion of a sphere such as in Kidambi \& Newton \cite{kidnew2}, is in
general not LA-integrable, careful numerical studies such as that of Knio
\emph{et al}. \cite{kncj} provide compelling evidence of the existence of
chaos for the three vortex problem in the half-plane.

From a perturbation perspective, we can obviously write the Hamiltonian
function in the form
\begin{equation}
H_{A}=H_{0}+H_{A,1}, \label{e25}%
\end{equation}
where
\begin{equation}
H_{A,1}:=%
{\sum\limits_{1\leq j\leq k\leq3}}
\kappa_{j}\kappa_{k}\log\left\vert \bar{z}_{j}-z_{k}\right\vert . \label{e26}%
\end{equation}
Now it is obvious from the nature of these various functions that there exist
extensive regions in $\mathbb{H}_{\#}$ at any positive distance away from the
boundary (the $x$-axis) where we have both $\left\vert H_{a}\right\vert
\ll\left\vert H_{0}\right\vert $ and $\left\vert \partial_{z}H_{a}\right\vert
\ll\left\vert \partial_{z}H_{0}\right\vert $, in keeping with our point of
view of treating (21) as a (small) perturbation of three vortex dynamics.

\medskip

\noindent\textsf{\textbf{B. Restricted four vortex dynamics}}

\smallskip

By the\emph{ restricted four vortex problem} (\emph{dynamics}), we mean the
motion of four point vortices in an ideal fluid in the (complex) plane
$\mathbb{C}$ $(=\mathbb{R}^{2})$, where three of the vortices, located at
points $z_{1}$, $z_{2}$ and $z_{3}$, have respective nonzero strengths
$\Gamma_{1}$, $\Gamma_{2}$ and $\Gamma_{3}$, while the fourth vortex at the
point $z_{4}$ has strength $\Gamma_{4}$ satisfying $\left\vert \Gamma
_{4}\right\vert \ll\mathfrak{m}:=\min\{\left\vert \Gamma_{1}\right\vert
,\left\vert \Gamma_{2}\right\vert ,\left\vert \Gamma_{3}\right\vert \}$ and we
have the freedom of making $\left\vert \Gamma_{4}\right\vert /\mathfrak{m}$ as
small as we wish. It follows from (4) that the complex Hamiltonian governing
equations are%
\begin{equation}
\kappa_{k}\dot{\bar{z}}_{k}=i%
{\sum\limits_{j=1,j\neq k}^{4}}
\kappa_{j}\kappa_{k}\left(  z_{j}-z_{k}\right)  ^{-1}=2i\partial_{z_{k}}%
H_{B},\qquad(1\leq k\leq4) \label{e27}%
\end{equation}
where the Hamiltonian function is%
\begin{equation}
H_{B}:=-%
{\sum\limits_{1\leq j<k\leq4}}
\kappa_{j}\kappa_{k}\log\left\vert z_{j}-z_{k}\right\vert . \label{e28}%
\end{equation}
Of course this can be put into the following real Hamiltonian form using the
same Poisson brackets as in (4):%
\begin{equation}
\dot{x}_{k}=\kappa_{k}^{-1}\partial_{y_{k}}H_{B}=\left\{
H_{B},x_{k}\right\}
,\quad\dot{y}_{k}=-\kappa_{k}^{-1}\partial_{x_{k}}H_{B}=\left\{  H_{B}%
,y_{k}\right\}  ,\qquad(1\leq k\leq4) \label{e29}%
\end{equation}
where analogous adjustments in the domain are obviously required to insure
smoothness. Just as in the case of (4), symmetry considerations show that (26)
has the motion invariants%
\begin{equation}
H_{B},\;J:=%
{\sum\limits_{k=1}^{4}}
\kappa_{k}z_{k},\;K:=%
{\sum\limits_{k=1}^{4}}
\kappa_{k}\left\vert z_{k}\right\vert ^{2}, \label{e30}%
\end{equation}
with the following independent integrals in involution
\begin{equation}
H_{B},\;K,\;L:=\left(  \mathfrak{R}J\right)  ^{2}+\left(  \mathfrak{I}%
J\right)  ^{2}. \label{e31}%
\end{equation}
These three independent integrals in involution represent the maximal number,
as it has been proved that the system (29) is not integrable in general (see
Bagrets \& Bagrets \cite{bb} and Ziglin \cite{ziglin}). We note, however, that
there are certain special cases of four vortex motion on a plane or sphere
that are LA-integrable, as shown in Aref \& Stremler \cite{arefstrem}, Borisov
\emph{et al}. \cite{bormamkil2}, Eckhardt \cite{eck}, Sakajo \cite{sakajo},
and Sokolovskiy \& Verron \cite{sokver}.

Our focus here is to treat (27) or (29) as a perturbation of the three vortex
problem. To be mathematically precise, we actually have to consider the system
as a perturbation of the three vortex problem embedded in $\mathbb{C}^{4}$, by
which we mean%
\begin{align}
\kappa_{k}\dot{\bar{z}}_{k}  &  =i%
{\sum\limits_{j=1,j\neq k}^{3}}
\kappa_{j}\kappa_{k}\left(  z_{j}-z_{k}\right)  ^{-1}=2i\partial_{z_{k}}%
H_{0},\qquad(1\leq k\leq3)\nonumber\\
\dot{\bar{z}}_{4}  &  =0=2i\partial_{z_{4}}H_{0}. \label{e32}%
\end{align}
Obviously this embedded system is also LA-integrable, and its dynamics
consists of a two-parameter infinity of copies of three vortex dynamics - one
for each constant value of $z_{4}$. This having been said, there is no harm in
the slight abuse of notation that we shall use from now on of treating the
restricted four vortex problem as a perturbation of the three vortex problem.
To highlight the perturbation aspect, we write the Hamiltonian function in the
form%
\begin{equation}
H_{B}=H_{0}+H_{B,1}, \label{e33}%
\end{equation}
where%
\begin{equation}
H_{B,1}:=-\sigma%
{\sum\limits_{j=1}^{3}}
\kappa_{j}\log\left\vert z_{j}-z_{4}\right\vert , \label{e34}%
\end{equation}
and we have replaced $\kappa_{4}$ by $\sigma$ to emphasize the fact that this
parameter may be chosen to be as small as necessary to insure the desired
perturbation behavior.

\medskip

\noindent\textsf{\textbf{C. Three slender coaxial vortex ring
dynamics}}

\smallskip

Consider three (circular) coaxial vortex rings of respective nonzero strengths
$\Gamma_{1}$, $\Gamma_{2}$ and $\Gamma_{3}$, with axis of symmetry the
$y$-axis, in motion in an ideal fluid in $\mathbb{R}^{3}$. The rings intersect
any meridian half-plane bounded along the $y$-axis in unique points
$(r_{1},y_{1})$, $(r_{2},y_{2})$ and $(r_{3},y_{3})$, respectively, where $r$
is the distance from the $y$-axis, and the motion of the rings is completely
determined by these points owing to the axial symmetry of the equations of
motion. It is convenient to set $z:=r^{2}+iy=s+iy$, so that the motion of the
three intersection points can be considered to be in the half-plane
$\mathbb{H}$. The complex Hamiltonian form of the dynamical equation of motion
of these points - de-singularized in a standard way to eliminate the usual
infinity in the self-induced velocity of the rings -can be expressed as
(\emph{cf}. Blackmore \& Knio \cite{blkn, blkn2}, Blackmore et al. \cite{bcw},
and Lamb \cite{lam})%
\begin{align}
\kappa_{k}\dot{\bar{z}}_{k}  &  =%
{\sum\limits_{j=1,j\neq k}^{3}} \kappa_{j}\kappa_{k}r_{j}r_{k}\left(
z_{j}-z_{k}\right) {\int\nolimits_{0}^{\pi/2}}
\Delta_{jk}^{-3/2}\cos2\alpha
d\alpha-i\kappa_{k}^{2}r_{k}^{-1}\left[
\log\left(  8r_{k}\delta_{c}^{-1}\right)  -0.558\right] \nonumber\\
&  -2i%
{\sum\limits_{j=1,j\neq k}^{3}}
\kappa_{j}\kappa_{k}r_{j}%
{\int\nolimits_{0}^{\pi/2}}
\left(  r_{j}-r_{k}\cos2\alpha\right)  \Delta_{jk}^{-3/2}d\alpha\label{e35}\\
&  =2i\partial_{z_{k}}H_{C},\quad(1\leq k\leq3)\nonumber
\end{align}
where $\delta_{c}$ is a very small positive number representing the common
core radii (used in the de-singularization procedure) of the vortex rings,
\begin{equation}
\Delta_{jk}:=\left(  r_{k}-r_{j}\right)  ^{2}+\left(  y_{k}-y_{j}\right)
^{2}+4r_{j}r_{k}\sin^{2}\alpha, \label{e36}%
\end{equation}
and the Hamiltonian function is%
\begin{equation}
H_{C}:=-2%
{\sum\limits_{j=1}^{3}} \kappa_{j}^{2}r_{j}\left[ \log\left(
8r_{j}\delta_{c}^{-1}\right)
-1.558\right]  -4%
{\sum\limits_{1\leq j<k\leq3}}
\kappa_{j}\kappa_{k}r_{j}r_{k}%
{\int\nolimits_{0}^{\pi/2}}
\Delta_{jk}^{-1/2}\cos2\alpha d\alpha. \label{e37}%
\end{equation}

This system can, employing the same Poisson bracket used in (4), be recast in
the real Hamiltonian form%
\begin{equation}
\dot{s}_{k}=\kappa_{k}^{-1}\partial_{y_{k}}H_{C}=\left\{
H_{C},s_{k}\right\}
,\quad\dot{y}_{k}=-\kappa_{k}^{-1}\partial_{s_{k}}H_{C}=\left\{  H_{C}%
,y_{k}\right\}  .\qquad(1\leq k\leq3) \label{e38}%
\end{equation}
It is easy to show that (37) has the two following independent integrals in
involution%
\begin{equation}
H_{C},\;G:=%
{\sum\limits_{j=1}^{3}}
\kappa_{j}s_{j}=%
{\sum\limits_{j=1}^{3}}
\kappa_{j}r_{j}^{2}. \label{e39}%
\end{equation}
However, there are no additional independent invariants in involution, as
proved in Bagrets \& Bagrets \cite{bb}, so although the dynamics of two
coaxial rings is LA-integrable, this is not true in general for the system
(35) or (38).

The formulation of (35) as a perturbation of three vortex dynamics is rather
more subtle than that of the two preceding examples. Formulas obtained by
Callegari \& Ting \cite{calting} indicate that, relative to a coordinate
system moving along the axis of symmetry with the overall translation velocity
of the ring configuration, the equations of motion of the rings are closely
approximated (modulo a constant factor) by those of three point vortices at
the intersection points of the rings in a meridian plane, when the rings are
sufficiently close to one another compared to their distance from the axis of
symmetry. More precisely, Blackmore \& Knio \cite{blkn} showed that in a
coordinate system moving with the \emph{center of vorticity }
\begin{equation}
z_{cv}:=\frac{\kappa_{1}z_{1}+\kappa_{2}z_{2}+\kappa_{3}z_{3}}{\kappa
_{1}+\kappa_{2}+\kappa_{3}}=J/\mathfrak{K}^{(1)}, \label{e40}%
\end{equation}
which is defined for any planar configuration of point vortices as long as
$\mathfrak{K}^{(1)}\neq0$, the equations of motion are Hamiltonian, with a
Hamiltonian function of the form%
\begin{equation}
\tilde{H}_{C}=-\left(  \mathfrak{R}z_{cv}/2\right)  \left[  H_{0}%
+H_{1}\right]  , \label{e41}%
\end{equation}
where $\left\vert H_{1}\right\vert /\left\vert H_{0}\right\vert
=o\left( 1/\log\rho\right)$ and $\left\vert
\partial_{z}H_{1}\right\vert /\left\vert
\partial_{z}H_{0}\right\vert =o\left(  \rho\right)$ as $\rho\rightarrow0$,
where $\rho$ is the ratio of the diameter of the configuration of meridian
plane points of intersection of the rings to the distance of $z_{cv}$ from the
axis of symmetry. As $\mathfrak{R}z_{cv}$ is a constant of motion of (35), the
above analysis establishes the three coaxial ring problem as a perturbation of
the three vortex problem for rings. We note that the above formulas require
that $\mathfrak{K}^{(1)}\neq0$, and, to simplify matters, we shall usually
assume hereafter that this is the case. One can prove the results that we
obtain in the sequel without this assumption, but certain rather
straightforward modifications are required in our methods of proof in some instances.

\bigskip

\noindent\textsf{\textbf{IV. QUALITATIVE REGULARITY RESULTS}}

\medskip

To get an idea of the type of results that we shall present in this section,
we first state a theorem for the examples of the preceding section for the
cases in which all of the vortex strengths have the same sign. The proof of
this theorem is either contained in, or follows directly from, the results in
Blackmore \& Knio \cite{blkn, blkn2}, Blackmore \emph{et al}. \cite{bcw},
Blackmore \emph{et al}. \cite{btk}, Blackmore \& Champanerkar \cite{blch},
Khanin \cite{khan}, and Ting \& Blackmore \cite{tb}.

\medskip

\noindent\textbf{Theorem 1.}\textit{ Suppose that in each of the perturbation
examples in the preceding section the parameters }$\kappa_{j}$\textit{ all
have the same sign, so that }$K^{(1)}\neq0$\textit{ and the center of
vorticity is defined. Then the following properties hold with respect to a
moving coordinate system with origin at the center of vorticity:}

\begin{itemize}
\item[(i)] \textit{For the three vortex problem in the half-plane (23)
governed by the Hamiltonian function }$H_{A}$\textit{ there is a set
of initial configurations of positive (Lebesgue) measure with the
diameter of the
configuration sufficiently small with respect to the distance from the }%
$x$\textit{-axis such that the motion is quasiperiodic on invariant KAM tori,
interspersed with periodic orbits.}

\item[(ii)] \textit{The restricted four vortex dynamics governed by }$H_{B}%
$\textit{ according to (29) exhibits, for a set of initial conditions of
positive measure, quasiperiodicity on invariant KAM tori along with periodic
orbits. This holds when the diameter of the initial configuration is
sufficiently small and }$\sigma$\textit{ is sufficiently small, or the fourth
vortex is sufficiently far from the three larger vortices, or for some
combination of both the conditions on }$\sigma$\textit{ and the distance of
the smaller vortex from the larger vortices.}

\item[(iii)] \textit{Three coaxial vortex ring dynamics generated by }$H_{C}%
$\textit{ in accordance with (38) exhibits the following behavior: There is a
set of positive measure of initial configurations with the distances among
them sufficiently small compared to their minimum distance to the axis of
symmetry that produces quasiperiodic motion on invariant KAM tori and periodic
orbits.}
\end{itemize}

\smallskip

\noindent\textit{Proof}. The proof, as indicated above, follows directly from
the cited papers. $\square$

\medskip

In the remainder of this section, we shall find rather general conditions on a
perturbation term $H_{1}$ of a Hamiltonian system%
\begin{equation}
\dot{\bar{z}}_{k}=\left\{  H,z_{k}\right\}  , \label{e42}%
\end{equation}
where $1\leq k\leq3$, or $1\leq k\leq4$ in the embedded version described
above for the restricted four vortex problem. Here the Poisson bracket is the
same as used in all of the previous equations, the Hamiltonian function is of
the form%
\begin{equation}
H=H_{0}+H_{1}, \label{e43}%
\end{equation}
and the assumptions are general enough to subsume the persistence of
regularity results of Theorem 1, as well as including the same types of
perturbations when the vortex strengths associated to the three vortex
dynamics generated by $H_{0}$ are allowed to differ in sign.

To accommodate all of the perturbations discussed in Section \textsf{III}, we
shall consider $H$ to be defined and analytic on an open subset of an
appropriate complex unitary space, such as $\mathbb{H}_{\#}^{3}$ for the
half-plane and coaxial ring problems, and $\mathbb{C}_{\#}^{4}$ for the
restricted four vortex problem. Moreover, we shall also include the
possibility of the perturbation $H_{1}$ depending on a real (freely chosen)
parameter $\sigma$ such that
\begin{equation}
\left\vert H_{1}/H_{0}\right\vert ,\left\vert \partial_{\mathbf{z}}%
H_{1}\right\vert /\left\vert \partial_{\mathbf{z}}H_{0}\right\vert
\rightarrow0 \label{e44}%
\end{equation}
uniformly as $\left\vert \sigma\right\vert \rightarrow0$ on any compact subset
of the domain on which these expressions are defined. It is convenient to
introduce some additional notation. \ We shall refer to a subset of the domain
of definition of (42) as \emph{ample} if it is of positive (Lebesgue) measure
in the domain or a submanifold of the domain of complex dimension at least
two. Moreover, we shall say that the system (or $H$) has the \emph{uniform
perturbation property} (\emph{UPP}) if there exists a sequence $\mathcal{S}%
=\{S_{n}:n\in\mathbb{N}\}$ of subsets of the domain of the system with
nonempty interiors such that $S_{n+1}\subset S_{n}$,
\begin{equation}
\left\vert H_{1}/H_{0}\right\vert \leq\frac{1}{\log(n+1)},\text{ and
}\left\vert \partial_{\mathbf{z}}H_{1}\right\vert /\left\vert \partial
_{\mathbf{z}}H_{0}\right\vert \leq\frac{1}{n+1} \label{eq45}%
\end{equation}
for every\textbf{ }$\mathbf{z}\in S_{n}$ for all $n\in\mathbb{N}$, where
$\mathbb{N}$ denotes the positive integers. In this case, we call the sequence
$\mathcal{S}$ a \emph{uniform perturbation filtration }(\emph{UPF}) for (42).
We are now in a position to state the main result of this section in a very
concise way.

\medskip

\noindent\textbf{Theorem 2.} \textit{Let system (42) satisfy the UPP as
described in (45), and also (44) when the perturbation depends on a parameter
as in the case of the restricted four vortex problem. Furthermore, suppose
that for some UPF, }$\mathcal{S}=\{S_{n}:n\in\mathbb{N}\}$\textit{, the system
satisfies the property that there is a companion sequence of compact subsets
of positive measure }$\mathcal{K}=\{K_{n}:n\in\mathbb{N}\}$ \textit{with
}$K_{n}\subset S_{n}$\textit{ for all }$n\in\mathbb{N}$\textit{ such that}
\begin{equation}
\mathbf{z}(0)\in K_{n+1}\Longrightarrow\mathbf{z}(t)\in S_{n} \label{e46}%
\end{equation}
\textit{for all }$t\geq0$\textit{ and for each sufficiently large }%
$n$\textit{, where }$\mathbf{z}(t)$\textit{ is the solution of (42)
initially at }$\mathbf{z}(0)$\textit{. Then there is an ample set of
initial conditions for the system (42), and a set of sufficiently
small values of }$\left\vert \sigma\right\vert $,\textit{ when there
is such a free parameter dependence as in the case of the restricted
four vortex problem, for which the dynamics includes quasiperiodic
motion on an ample set of invariant KAM tori interspersed with a
countable collection of periodic orbits.}

\smallskip

\noindent\textit{Proof}. Our argument relies heavily on what might be called
the limit KAM theorem and a useful generalization of the Poincar\'{e}-Birkhoff
theorem, which were successfully employed by Blackmore and his collaborators
\cite{blkn, bcw, blch} to prove the existence of, respectively, an ample set
of invariant tori and periodic orbits for several examples of vortex dynamics
problems of the type under consideration.

First, following the same approach as in those papers, we rewrite the
Hamiltonian in terms of the action-angle coordinates associated to the
unperturbed LA-integrable system as%
\begin{equation}
\mathcal{H}\left(  \Lambda,\Theta\right)  =\mathcal{H}_{0}\left(
\Lambda\right)  +\mathcal{H}_{1}\left(  \Lambda,\Theta\right)  , \label{e47}%
\end{equation}
where the action and angle vectors are, respectively, $\Lambda=(\Lambda
_{1},...,\Lambda_{k})$ and $\Theta=$ $(\Theta_{1},...,\Theta_{k})$, with
$2\leq k\leq4$. We have used slightly different notation for the various
Hamiltonians to underscore the fact that it may be necessary, as in the case
of coaxial vortex rings, to reduce the number of degrees of freedom by one in
order to insure the nondegeneracy of $\mathcal{H}_{0}$, which naturally alters
the forms of the original terms of the Hamiltonian (43). This also explains
why we have indicated that the number of degrees of freedom $k$ in (47) can
assume values between two and four. More specifically, $k=2$ for the coaxial
vortex ring problem, $k=3$ for three vortex in a half-plane dynamics, and
$k=4$ for the restricted four vortex problem (in which case we recall that the
three vortex problem must be embedded in $\mathbb{C}^{4}$).

For our purposes it is not necessary to know the exact form of $\mathcal{H}%
_{0}\left(  \Lambda\right)  $, which is very difficult to deduce in general;
we need only take note of the following readily verifiable properties that
hold in any compact subset of the domain in which (42) is analytic:%
\begin{equation}
\Delta_{1}:=\det\left(  \frac{\partial^{2}\mathcal{H}_{0}\left(
\Lambda\right)  }{\partial\Lambda_{i}\partial\Lambda_{j}}\right)  \neq0,
\label{e48}%
\end{equation}%
\begin{equation}
\Delta_{2}:=\det\left(
\begin{array}
[c]{cc}%
\frac{\partial^{2}\mathcal{H}_{0}\left(  \Lambda\right)  }{\partial\Lambda
_{i}\partial\Lambda_{j}} & \frac{\partial\mathcal{H}_{0}\left(  \Lambda
\right)  }{\partial\Lambda_{i}}\\
\frac{\partial\mathcal{H}_{0}\left(  \Lambda\right)  }{\partial\Lambda_{j}} &
0
\end{array}
\right)  \neq0, \label{e49}%
\end{equation}
and%
\begin{equation}
\mathcal{H}_{0}\left(  \Lambda\right)  /\Delta_{1},\Delta_{1}/\Delta
_{2}\rightarrow0 \label{e50}%
\end{equation}
as $\left\vert \Lambda\right\vert ^{2}:=\Lambda_{1}^{2}+\cdot\cdot
\cdot+\Lambda_{k}^{2}\rightarrow0$ (\emph{cf}. \cite{blkn}).

Translating the hypotheses to the transformed Hamiltonian system generated by
$\mathcal{H}\left(  \Lambda,\Theta\right)  $; namely%
\begin{equation}
\dot{\Lambda}_{j}=\partial_{\Theta_{j}}\mathcal{H}=\partial_{\Theta_{j}%
}\mathcal{H}_{1},\;\dot{\Theta}_{j}=-\partial_{\Lambda_{j}}\mathcal{H}%
=-\partial_{\Lambda_{j}}\mathcal{H}_{0}-\partial_{\Lambda_{j}}\mathcal{H}_{1},
\label{e51}%
\end{equation}
we infer the existence for any $\epsilon>0$ of a compact set $K_{\epsilon}$ of
positive measure of initial conditions for (51) such that%
\begin{equation}
\left\vert \mathcal{H}_{1}\left(  \Lambda(t),\Theta(t)\right)  /\mathcal{H}%
_{0}\left(  \Lambda(t)\right)  \right\vert ,\left\vert \partial_{\left(
\Lambda,\Theta\right)  }\mathcal{H}_{1}\left(  \Lambda(t),\Theta(t)\right)
\right\vert /\left\vert \partial_{\Lambda}\mathcal{H}_{0}\left(
\Lambda(t)\right)  \right\vert \leq\epsilon\label{e52}%
\end{equation}
for all $t\geq0$. Consequently, owing to the nondegeneracy condition (48),
together with (50) in the limiting case for $\left\vert \Lambda\right\vert
\rightarrow0$ , it follows from the KAM theorem that there exists an ample set
of invariant (real) $k$-dimensional tori for (42). If $k=2$, the isoenergetic
nondegeneracy condition (49), coupled with (50) when $\left\vert
\Lambda\right\vert \rightarrow0$, is enough to insure the existence of
periodic orbits. For $k>2$, the existence of cycles for (42) can be proven,
for example, by using the generalization of the Poincar\'{e}-Birkhoff fixed
point theorem employed by Blackmore \emph{et al}. \cite{bcw}. Thus, the proof
is complete. $\square$

\medskip

It is a relatively straightforward task to demonstrate that the three types of
perturbations described in the preceding section satisfy the above hypotheses,
thereby demonstrating the existence of ample regular dynamics regimes (like
those for LA-integrable systems) in each of those cases; to wit, the following
result is a direct consequence of Theorem 2.

\medskip

\noindent\textbf{Corollary 1. }\textit{The three vortex in a
half-plane system (21), which we denote as (A), the restricted four
vortex problem (27), identified as (B), and the three coaxial vortex
ring system (35), denoted as (C), all satisfy the hypotheses of
Theorem 2, so each of these examples exhibits ample sets of initial
configurations leading to quasiperiodic flows on invariant KAM tori
accompanied by periodic orbits. More specifically, admissible UPF's
of initial conditions leading to regular motion for each of these
examples can be described as follows for the elliptic (ET),
hyperbolic (HT), and parabolic (PT) types of the unperturbed system:
}

\begin{itemize}
\item[(ET)] \textit{For (A) there are four kinds of UPF's , denoted as }%
$\mathcal{S}_{A}^{(0)}$\textit{ , }$\mathcal{S}_{A}^{(3)},\mathcal{S}_{A}^{(2)}$\textit{ and }$\mathcal{S}_{A}^{(1)}%
$\textit{, with the following characterizations:
}$\mathcal{S}_{A}^{(0)}$\textit{ is comprised of approximately
equilateral configurations (associated to the point }$E$\textit{ or
}$E^{\ast}$\textit{ described in Section II) of the three vortices
with diameter sufficiently small compared with their distance from
the }$x$\textit{-axis; }$\mathcal{S}_{A}^{(3)}$\textit{ consists of
configurations with }$\left\vert z_{1}-z_{2}\right\vert
\ll\left\vert z_{3}-z_{2}\right\vert
\simeq\left\vert z_{3}-z_{1}\right\vert $\textit{ (associated with }$Q_{3}%
$\textit{ defined in Section II) of diameter sufficiently small compared with
the distance of the configuration from the }$x$\textit{-axis; }$\mathcal{S}_{A}^{(2)}%
$\textit{, associated with }$Q_{2}$\textit{, is characterized by
making the
obvious changes in the definition of }$\ \mathcal{S}_{A}^{(3)}$\textit{; and }%
$\mathcal{S}_{A}^{(1)}$\textit{, corresponding to }$Q_{1}$\textit{,
is defined by making the evident revisions of the description of
}$\mathcal{S}_{A}^{(3)}$\textit{. There are
also four kinds of UPF's for (B), }$\mathcal{S}_{B}^{(0)}$\textit{, }$\mathcal{S}_{B}^{(3)}%
$\textit{, }$\mathcal{S}_{B}^{(2)}$\textit{ and
}$\mathcal{S}_{B}^{(1)}$\textit{, defined analogously to
}$\mathcal{S}_{B}^{(\ast)}$\textit{ ,
}$\mathcal{S}_{A}^{(3)},\mathcal{S}_{A}^{(2)}$\textit{ and
}$\mathcal{S}_{A}^{(1)}$\textit{, respectively, wherein the ratio of
the diameter of
the configuration of the main three vortices to the distance from the }%
$x$\textit{-axis is replaced by the ratio of the diameter to the
distance of the configuration from the fourth small vortex; and four
kinds for (C), }$\mathcal{S}_{C}^{(\ast)}$\textit{,
}$\mathcal{S}_{C}^{(3)}$\textit{, }$\mathcal{S}_{C}^{(2)}$\textit{
and }$\mathcal{S}_{C}^{(1)}$\textit{, also defined analogously to
}$\mathcal{S}_{A}^{(\ast)}$\textit{ ,
}$\mathcal{S}_{A}^{(3)},\mathcal{S}_{A}^{(2)}$\textit{ and
}$\mathcal{S}_{A}^{(1)}$\textit{, respectively, where the point
vortices are replaced by the intersection points of the rings with
any given meridian plane, and the distance from the axis of symmetry
replaces the distance from the }$x$\textit{-axis.}

\item[(HT)] \textit{There is at least one kind of UPF for (A) corresponding to
}$\mathcal{S}_{A}^{(3)}$\textit{ as described above, and as many as
four others associated with }$\mathcal{S}_{A}^{(2)}$\textit{ and
}$\mathcal{S}_{A}^{(1)}$\textit{, and
possibly two others corresponding to the points }$Q_{4}$\textit{ and }$Q_{5}%
$\textit{ when they are centers, as discussed in Section II. For (B)
and (C)
there are analogous UPF's that always include those associated to }%
$\mathcal{S}_{B}^{(3)}$\textit{ and
}$\mathcal{S}_{C}^{(3)}$\textit{, respectively.}

\item[(PT)] \textit{For (A), (B) and (C) there is always one UPF corresponding
to }$\mathcal{S}_{A}^{(3)}$\textit{,
}$\mathcal{S}_{B}^{(3)}$\textit{ and
}$\mathcal{S}_{C}^{(3)}$\textit{, respectively, and possibly
additional UPF associated to any other possible centers as described
in the trilinear phase plane characterization of the dynamics
covered in Section II.}
\end{itemize}

\smallskip

\noindent\textit{Proof}. We shall provide detailed arguments only for the
three vortex system in the half-plane, either when the strengths of all
vortices have the same sign (the elliptic type), or they differ in sign (which
includes both the hyperbolic and parabolic types). The verifications for the
restricted four vortex, and three coaxial vortex ring problems for the various
types of the unperturbed (three vortex) system can be obtained analogously by
straightforward - but rather lengthy - calculations, and we shall leave the
details to the reader, noting that the desired results when the unperturbed
system is of elliptic type actually follow directly from Theorem 1.

Even though, as mentioned above, the conclusions we are seeking are a direct
consequence of Theorem 1 when the unperturbed system is of elliptic type, we
shall include a proof here, since it will be helpful in pointing the way to
establishing sufficient conditions for the hyperbolic and parabolic types.
Recalling that there is no loss of generality in assuming that $\kappa_{1}%
\geq\kappa_{2}\geq\kappa_{3}>0$ when the unperturbed system is elliptic, we
first rewrite the constants of motion (24) of (23) in a form more useful to
our purposes; namely%
\begin{align}
H_{A}  &  =c_{1}=\kappa_{1}^{2}\log\left(  2y_{1}\right)  +\kappa_{2}^{2}%
\log\left(  2y_{2}\right)  +\kappa_{3}^{2}\log\left(  2y_{3}\right)
-\kappa_{1}\kappa_{2}\log\left[  \frac{(x_{1}-x_{2})^{2}+(y_{1}-y_{2})^{2}%
}{(x_{1}-x_{2})^{2}+(y_{1}+y_{2})^{2}}\right]  ^{1/2}-\nonumber\\
&  \kappa_{1}\kappa_{3}\log\left[  \frac{(x_{1}-x_{3})^{2}+(y_{1}-y_{3})^{2}%
}{(x_{1}-x_{3})^{2}+(y_{1}+y_{3})^{2}}\right]  ^{1/2}-\kappa_{2}\kappa_{3}%
\log\left[  \frac{(x_{2}-x_{3})^{2}+(y_{2}-y_{3})^{2}}{(x_{2}-x_{3}%
)^{2}+(y_{2}+y_{3})^{2}}\right]  ^{1/2},\label{e53}\\
\mathfrak{I}J  &  =c_{2}=\kappa_{1}y_{1}+\kappa_{2}y_{2}+\kappa_{3}%
y_{3,}\nonumber
\end{align}
and define $M\left(  c_{1},c_{2}\right)  $ to be the set of points in
$\mathbb{H}_{\#}^{3}$ satisfying the pair of equations (53).

It is straightforward to show directly from the form of these defining
equations that given any arbitrarily large and small positive number,
respectively, $\lambda$ and $\nu$, there exist $\lambda_{\ast}=\lambda_{\ast
}(\lambda,\nu)\geq\lambda$ and $0<\nu_{\ast}=\nu_{\ast}(\lambda,\nu)<\nu$ such
that the following property is satisfied if $c_{1}$ and $c_{2}$ are chosen to
be sufficiently large positive numbers: let the initial positions of the
vortices, $z_{1}(0)=x_{1}(0)+iy_{1}(0)$, $z_{2}(0)=x_{2}(0)+iy_{2}(0)$,
$z_{3}(0)=x_{3}(0)+iy_{3}(0)$ be approximately in the shape of an equilateral
triangle (corresponding to $E$ in Fig.1) of diameter less than $\nu_{\ast}$,
with $y_{1}(0),$ $y_{2}(0),y_{3}(0)\geq\lambda_{\ast}$, and define $M_{\ast
}\left(  c_{1},c_{2}\right)  $ to be the component of $M\left(  c_{1}%
,c_{2}\right)  $ containing the initial point $\mathbf{z}(0)$. Then for all
configurations $(z_{1},z_{2},z_{3})$ in the component $M_{\ast}\left(
c_{1},c_{2}\right)  $, the diameter is less than or equal to $\nu$, and
$y_{1},$ $y_{2},y_{3}\geq\lambda$. Owing to the connectedness of orbits of
(23), it must therefore follow that the diameter of the configuration
$\mathbf{z}(t)=(z_{1}(t),z_{2}(t),z_{3}(t))$ is less than or equal to $\nu$,
and $y_{1}(t),$ $y_{2}(t),y_{3}(t)\geq\lambda$ for all $t\geq0$, where
$\mathbf{z}(t)$ is the solution of (23) with the specified initial condition.
Whence the construction of the desired UPF associated to these sets is a
simple matter.

If the unperturbed system is hyperbolic or parabolic, then as indicated in
Section \textsf{II}, we may - and do - assume without loss of generality that
$\kappa_{1}\geq\kappa_{2}>0>\kappa_{3}$. This puts a slightly different
complexion on the system of equations (53), which must be satisfied by all
solutions of (23). However, not so different that we cannot use the same kind
of argument as for the elliptic type system (given above) modulo a few rather
obvious modifications. In fact, we can take our cue for the necessary
adjustments by recalling our discussion of the trilinear phase portraits for
the three vortex problem in Section \textsf{II}. It is not difficult to see
from a close inspection of (53) for the case when $\kappa_{3}$ is negative -
in which we concentrate on those terms containing the negative vortex strength
- that we can simply change the initial condition on the vortex configuration
to conform to the center $Q_{3}$ (see Figs 1 and 2). More precisely, we need
only change the description of the initial configuration in the preceding
paragraph so that $\left\vert z_{1}(0)-z_{2}(0)\right\vert \ll$ $\left\vert
z_{1}(0)-z_{3}(0)\right\vert \simeq$ $\left\vert z_{2}(0)-z_{3}(0)\right\vert
$ in order to obtain orbits of (23) staying in a component of $M\left(
c_{1},c_{2}\right)  $ analogous to $M_{\ast}\left(  c_{1},c_{2}\right)  $.
Just as in the previous paragraph, this leads directly to the desired UPF's
for the hyperbolic and parabolic cases, and completes the proof for three
vortex dynamics in the half-plane. $\square$

\medskip

Before moving on to a study of perturbations of specific periodic orbits of
three vortex dynamics, we note that the following generalization of Theorem 2
can be easily proved by using essentially the same arguments as in its proof
given above.

\medskip

\noindent\textbf{Theorem 3.}\textit{ Suppose the system (42) is a perturbation
of a general LA-integrable Hamiltonian system generated by }$H_{0}$\textit{,
and that }$H$\textit{ is defined and analytic on an open subset of
}$\mathbb{C}^{m}$\textit{, with }$m\geq1$\textit{. In addition, assume that
(42) satisfies the hypotheses of Theorem 2, and that the usual KAM
nondegeneracy condition holds for }$H_{0}$\textit{ in an admissible UPF. Then
there exists an ample set of initial conditions (and small values of a
parameter, if pertinent) such that the system exhibits quasiperiodic flows on
an ample collection of invariant tori, together with periodic orbits.}

\textit{ \bigskip}

\noindent\textsf{\textbf{V. PERSISTENCE OF PERIODIC ORBITS}}

\medskip

In the preceding section we proved theorems demonstrating that various types
of, generally non-integrable, Hamiltonian perturbations of certain
LA-integrable systems - including that governing the motion of three point
vortices in an ideal fluid in the complex plane - have ample dynamical regimes
exhibiting the regularity properties that characterize integrable systems;
namely quasiperiodic motion on invariant tori and periodic orbits.

Such qualitative results as those already obtained, which are essentially
existence theorems, beg the question of how closely such behaviors of the
perturbed system approximate those of the unperturbed LA-integrable system
associated to $H_{0}$? This can be viewed as a more quantitative classical
perturbation theory related query about the qualitative entities whose
existence has been proven by more modern methods in symplectic dynamics. In
this section we provide a partial answer to this question as it relates to
periodic orbits, which is embodied in the following result.

\medskip

\noindent\textbf{Theorem 4.} \textit{Let the Hamiltonian system defined by
(42) and (43) be a perturbation of three vortex dynamics generated by }$H_{0}%
$\textit{ with respect to a coordinate system moving with the center of
vorticity. Suppose that the closed curve }$C$\textit{ represents a periodic
orbit of the unperturbed three vortex system with Hamiltonian function }%
$H_{0}$\textit{. Define the (compact) tubular neighborhood }$\mathfrak{T}%
_{\epsilon}(C)$\textit{ for each positive }$\epsilon$\textit{ as}%
\[
\mathfrak{T}_{\epsilon}(C):=\left\{  \mathbf{z}:\Delta(\mathbf{z}%
,C)\leq\epsilon\right\}  ,
\]
\textit{where }$\Delta$\textit{ denotes the usual (unitary) distance function.
Also define}%
\[
\lambda_{\epsilon}(C)=\max\left\{  \left\vert H_{1}/H_{0}\right\vert
,\left\vert \partial_{\mathbf{z}}H_{1}\right\vert /\left\vert \partial
_{\mathbf{z}}H_{0}\right\vert :\mathbf{z}\in\mathfrak{T}_{\epsilon
}(C)\right\}  .
\]

\textit{Then if for a given }$\epsilon$\textit{ small enough to
insure that }$\mathfrak{T}_{\epsilon}(C)$\textit{ is a smooth
(}$=C^{\infty}$\textit{) submanifold of the domain in which the
system (42) is smooth, the quantity }$\lambda_{\epsilon}(C)$\textit{
is sufficiently small, (42) has a periodic orbit
}$\tilde{C}$\textit{ in }$\mathfrak{T}_{\epsilon}(C).$

\smallskip

\noindent\textit{Proof}. First select a point $\mathbf{z}_{0}\in C$, and let
$\mathcal{E}_{0}$ and $\mathcal{E}$ be the energy hypersurface, respectively,
of the unperturbed system generated by $H_{0}$ and the perturbed system
generated by $H$, which contain the point $\mathbf{z}_{0}$. Both
$\mathcal{E}_{0}$ and $\mathcal{E}$ have a common odd real dimension, which we
denote as $2m-1$. Let $\Delta_{r}$ denote the Riemannian metric in $E$ induced
by the metric $\Delta$, and define%
\[
B_{\delta}:=\left\{  \mathbf{z}\in\mathcal{E}:\Delta_{r}\left(  \mathbf{z}%
,\mathbf{z}_{0}\right)  \leq\delta\right\}
\]
for small positive values of $\delta$. This is clearly a $(2m-1)$-ball in
$\mathcal{E}$ for all $\delta\leq\delta_{0}$ sufficiently small.

We let $\varphi_{t}$ represent the flow generated by (42). The following
properties follow directly from standard results on differential equations
(see e.g. \cite{kathas}) when $\lambda_{\epsilon}(C)$ is chosen to be
sufficiently small: There exist a transversal $\Sigma$ through $\mathbf{z}%
_{0}$ for the system (42) in $\mathcal{E}$ and $0<\delta_{2}<\delta_{1}$
$<\delta_{0}$ such that (i) each of the sets $\beta_{\delta}:=B_{\delta}%
\cap\Sigma$ is a $2(m-1)$-ball in $\mathcal{E}$ whenever
$0<\delta\leq \delta_{1}$, (ii) the orbits of (42) initially on
$\beta_{\delta_{2}}$ all pass through the interior of
$\beta_{\delta_{1}}$ in finite (positive time) generating a
Poincar\'{e} map $P:\beta_{\delta_{2}}\rightarrow$ $\beta
_{\delta_{1}}$ with $P\left(  \beta_{\delta_{2}}\right)  \subset
\mathrm{interior}(\beta_{\delta_{2}})$, (iii) the radial geodesic
curves in $E$ emanating from $\varphi_{t}(\mathbf{z}_{0})$ all
transversely intersect the boundary of
$\varphi_{t}(\beta_{\delta_{2}})$ in unique points for all $0\leq
t\leq t_{m}$, where $t_{m}$ is the maximum first return time to
$\Sigma$ among all those for the flow of $\beta_{\delta_{1}}$, and
(iv) the flow of $\beta_{\delta_{1}}$ generated by (42) through the
first return time to $\Sigma$ is contained in
$\mathfrak{T}_{\epsilon}(C)$.

Our proof is obviously complete if $P$ has a fixed point, so assume on the
contrary that this is not the case. As $\beta_{\delta_{1}}$ and $\beta
_{\delta_{2}}$ are both (real) odd-dimensional balls, it follows from
properties (i)-(iv) that we can modify (42) in the flow of $\beta_{\delta_{1}%
}$ $\setminus\beta_{\delta_{2}}$ in order to obtain an extension $\hat{P}$ of
$P$, such that $\hat{P}$ maps $\beta_{\delta_{1}}$ into itself and has no
fixed points in $\beta_{\delta_{1}}$ $\setminus\beta_{\delta_{2}}$. Thus, this
smooth self-mapping $\hat{P}$ of $\beta_{\delta_{1}}$ has no fixed points at
all, which contradicts the Brouwer fixed point theorem. Accordingly we
conclude that $P$ must itself have a fixed point, which corresponds to a cycle
$\tilde{C}$ of (42) in $\mathfrak{T}_{\epsilon}(C)$, so the proof is complete.
$\square$

\medskip

Considering the generality of the above result, and the relative simplicity of
the proof, it seems as though this should have certainly been discovered.
However, the authors were unable to find such a result in the literature,
although it appears that it could be derived rather directly from certain
results that prove the existence of periodic solutions for Hamiltonian systems
employing index theory (\emph{cf}. Blackmore \& Wang \cite{blwa}, Gol\'{e}
\cite{gole} and Josellis \cite{jos}). It is interesting to take note of the
very special case where $\partial_{\mathbf{z}}H_{1}$ in the above theorem
takes the form $\partial_{\mathbf{z}}H_{1}=\phi(\mathbf{z})\partial
_{\mathbf{z}}H_{0}$, where $\phi$ is a smooth, real valued function whose
magnitude can be made arbitrarily small in the tubular neighborhood
$\mathfrak{T}_{\epsilon}(C)$. Then as one can readily see by making the
obvious transformation of the time parameter, the cycle of the perturbed
system is not just an approximation of $C$, it is identical with it, and the
local flow of the perturbed system is identical with that of the unperturbed
system (modulo a change of parametrization).

Theorem 4 can be applied directly to the perturbations of three vortex
dynamics that we have been considering to identify conditions under which the
perturbation has a periodic orbit close to one for the three vortex problem.
We leave the very straightforward proof of the next result to the reader.

\medskip

\noindent\textbf{Corollary 2. }\textit{The hypothesis and
conclusions of Theorem 4 regarding the existence of a periodic orbit
of the perturbed system close to a periodic orbit
}$C:z=z(t)$\textit{, }$0\leq t\leq p$\textit{, of the unperturbed
(three vortex) system hold for (A) three vortices in a half-plane,
(B) restricted four vortex dynamics, and (C) three slender coaxial
ring dynamics obtain under the following conditions:}

\begin{itemize}
\item[(a)] \textit{Each of the coordinates of }$z(t)$\textit{ on }$C$\textit{
is sufficiently distant from the }$x$\textit{-axis in the complex plane.}

\item[(b)] \textit{For the coordinates }$z(t)=(z_{1}(t),z_{2}(t),z_{3}%
(t),z_{4}(t))\in C$\textit{, the distance in the complex plane }$\mathbb{C}%
$\textit{ from }$z_{4}(0)$\textit{ to the set }$\{z\in\mathbb{C}%
:z=z_{1}(t),z_{2}(t)$\textit{\ }$or\;\,z_{3}(t),0\leq t\leq p\}$\textit{ is
sufficiently large, or the magnitude of the strength of the fourth vortex is
sufficiently small, or a suitable combination of both of these conditions is
\ enforced.}

\item[(c)] \textit{The coordinates of }$z(t)$\textit{ on }$C$\textit{
representing the points of intersection of the rings with a meridian plane are
sufficiently far from the axis of symmetry (represented by the }%
$y$\textit{-axis in this plane).}
\end{itemize}

We note that the conditions given in Corollary 2 by no means exhaust all
possible situations where the perturbed system has a periodic orbit close to a
periodic orbit of the unperturbed system. For example, it is easy to see that
in the restricted four vortex problem there are many such cases where the
fourth small vortex is neither particularly small nor very distant from the
three large vortices: Simply consider any configuration of the large vortices
that generates a periodic solution of the three vortex problem, and place the
fourth vortex, of any strength, at the center of vorticity of the larger
vortices. Then the fourth vortex remains fixed, and the motion of the three
larger vortices is unaffected by its presence.

\bigskip

\noindent\textsf{\textbf{VI. SIMULATIONS OF DYNAMICS}}

\medskip

In this section we provide numerical examples that illustrate the
persistence of regular motion for different types of perturbations,
as well as breakdown of regularity as large perturbations are
considered and the hypotheses of our theoretical results are
accordingly violated. As mentioned earlier, we consider half-plane
dynamics, restricted four vortex perturbations, and slender coaxial
vortex ring dynamics , which we refer to below as type A, B, and C,
respectively.

Our simulations cover just a small sample of possible cases for the
various perturbations, and they are presented in two basic types of
graphical forms: As trajectories in the plane or half-plane, or as
Poincar\'{e} sections, which are especially well suited to
illuminating transitions from regular to chaotic motion. In
particular, for the unperturbed system we present plots of the
trajectories of the three vortices in the plane, juxtaposed with the
corresponding trilinear phase plane stationary point. If any of the
perturbation types satisfy the hypotheses of Corollaries 1 and 2, we
expect that plots of the trajectories of their vortex elements (with
respect to a coordinate system moving with the configuration) would
show a small variation of the plots for the unperturbed system. Also
included for type A perturbations are Poincar\'{e} maps, along with
the corresponding Poincar\'{e} map for the unperturbed system,
showing that there is a transition to chaotic motion as the array
starts closer to the boundary of the half-plane compared to the
diameter of the initial configuration of vortices. For type C
perturbations, we present Poincar\'{e} maps exhibiting strong
regularity when the slender coaxial rings are very close to one
another compared to the distance of the configuration from the axis
of symmetry, as expected in view of Corollary 1. For the restricted
four vortex problem, a pair of Poincar\'{e} maps shows how the
motion tends to be more chaotic as the initial position and strength
of the fourth small vortex, respectively, starts closer to the
initial group of three larger vortices and grows in comparison to
the strengths of this group. This behavior is entirely consistent
with Corollaries 1 and 2.

\begin{figure}[h]
\centerline{
\includegraphics[angle=0,height=63mm,draft=false]
{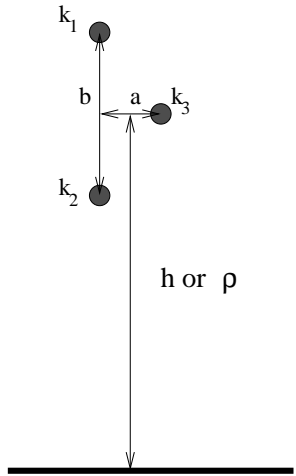}
}
\caption{Initial configuration for type I and II perturbations.}
\label{setup}
\end{figure}

For type A and C perturbations, we consider the setup illustrated in
Fig.~\ref{setup}.  As shown in the figure, the geometry is specified
in terms of two parameters, $b$ and $a$, while the strength of the
perturbation is reflected by the height $h$ for type A, or the mean
radius $\rho$ for type C.  Introducing the parameter $c \equiv
\sqrt{a^2+b^2}$, an analysis was conducted of the nine cases
summarized in Table~\ref{tab:inputs}.

\bigskip

\begin{table}[h]
\begin{center}
\begin{tabular}{ccccccc}
\hline
Case & $ b$ & $c$ & $\kappa_1$ & $\kappa_2$ & $\kappa_3$ & $\mathfrak{K}^{(2)}$ \\
\hline
1          &   0.1  &   1   &     2      &     1      &    0.5    &  3.5 \\
2          &   0.1  &   1   &     2      &     1      &  -2/3     &   0   \\
3          &   0.1  &   1   &     2      &     1      &  -0.8     & -0.4 \\
4          &   2/7  &   1   &     2      &     1      &    0.5    &  3.5 \\
5          &   2/7  &   1   &     2      &     1      &  -2/3     &   0   \\
6          &   2/7  &   1   &     2      &     1      &  -0.8     & -0.4 \\
7          &   0.49  &   1   &     2      &     1      &    0.5    &  3.5 \\
8          &   0.49  &   1   &     2      &     1      &  -2/3     &   0   \\
9          &   0.49  &   1   &     2      &     1      &  -0.8     & -0.4 \\
\hline
\end{tabular}
\caption{Summary of inputs.  The calculated value of
$\mathfrak{K}^{(2)}$ is also reported.} \label{tab:inputs}
\end{center}
\end{table}

\bigskip

\begin{figure}[phtb]
\centerline{
\includegraphics[angle=0,height=60mm,draft=false]
{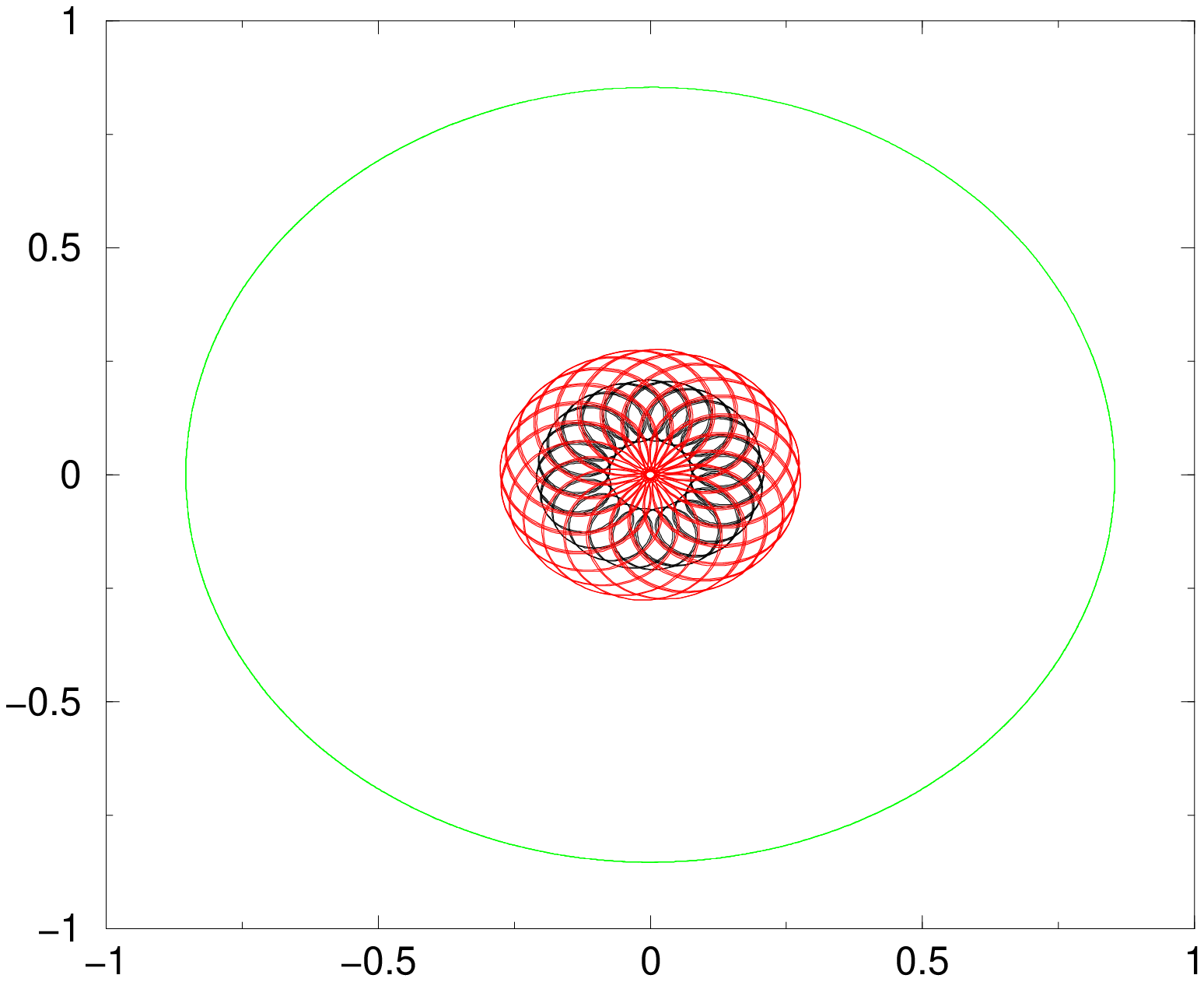} \hspace{5mm}
\includegraphics[angle=0,height=60mm,draft=false]
{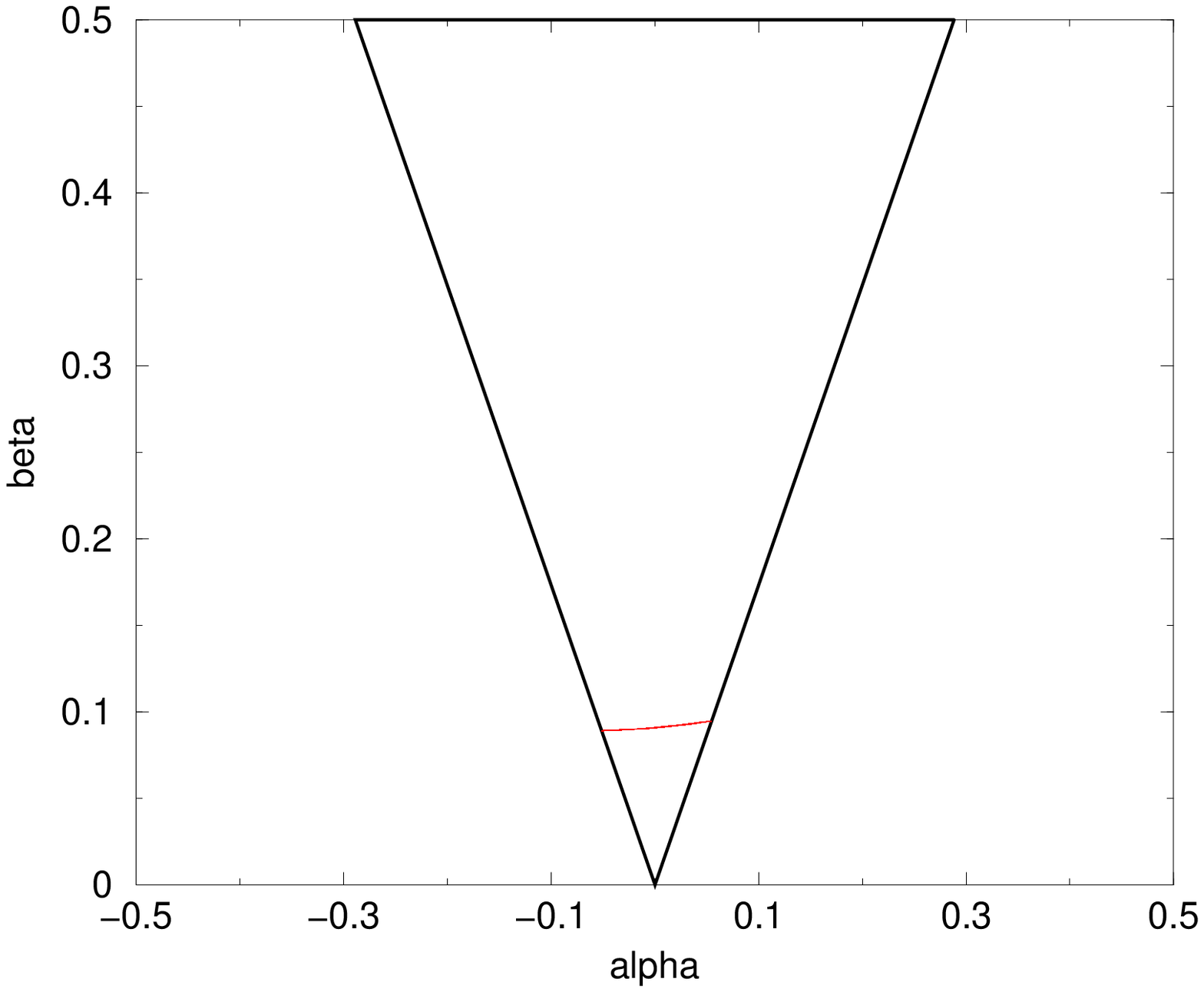}
}
\centerline{
\includegraphics[angle=0,height=60mm,draft=false]
{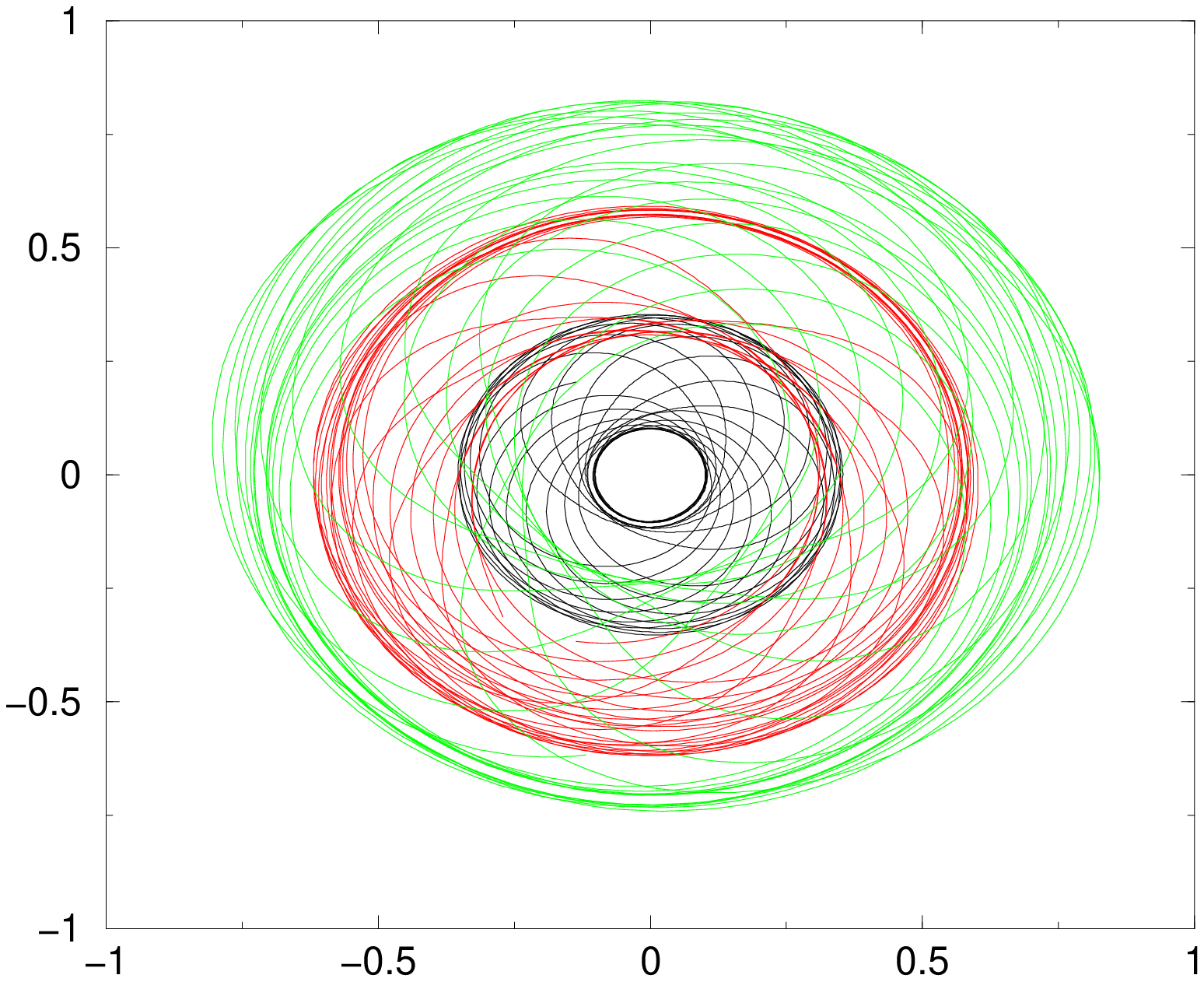} \hspace{5mm}
\includegraphics[angle=0,height=60mm,draft=false]
{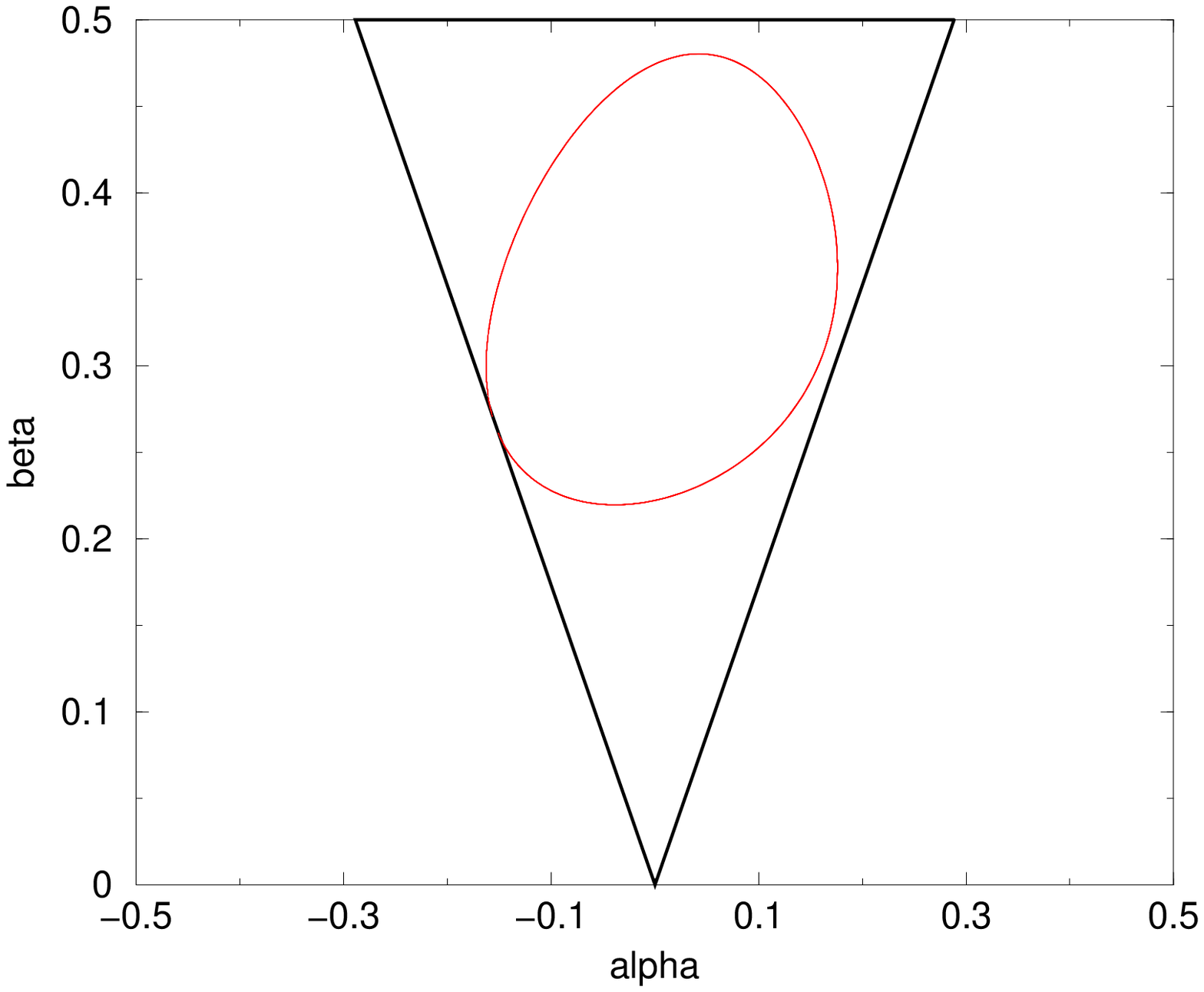}
}
\centerline{
\includegraphics[angle=0,height=60mm,draft=false]
{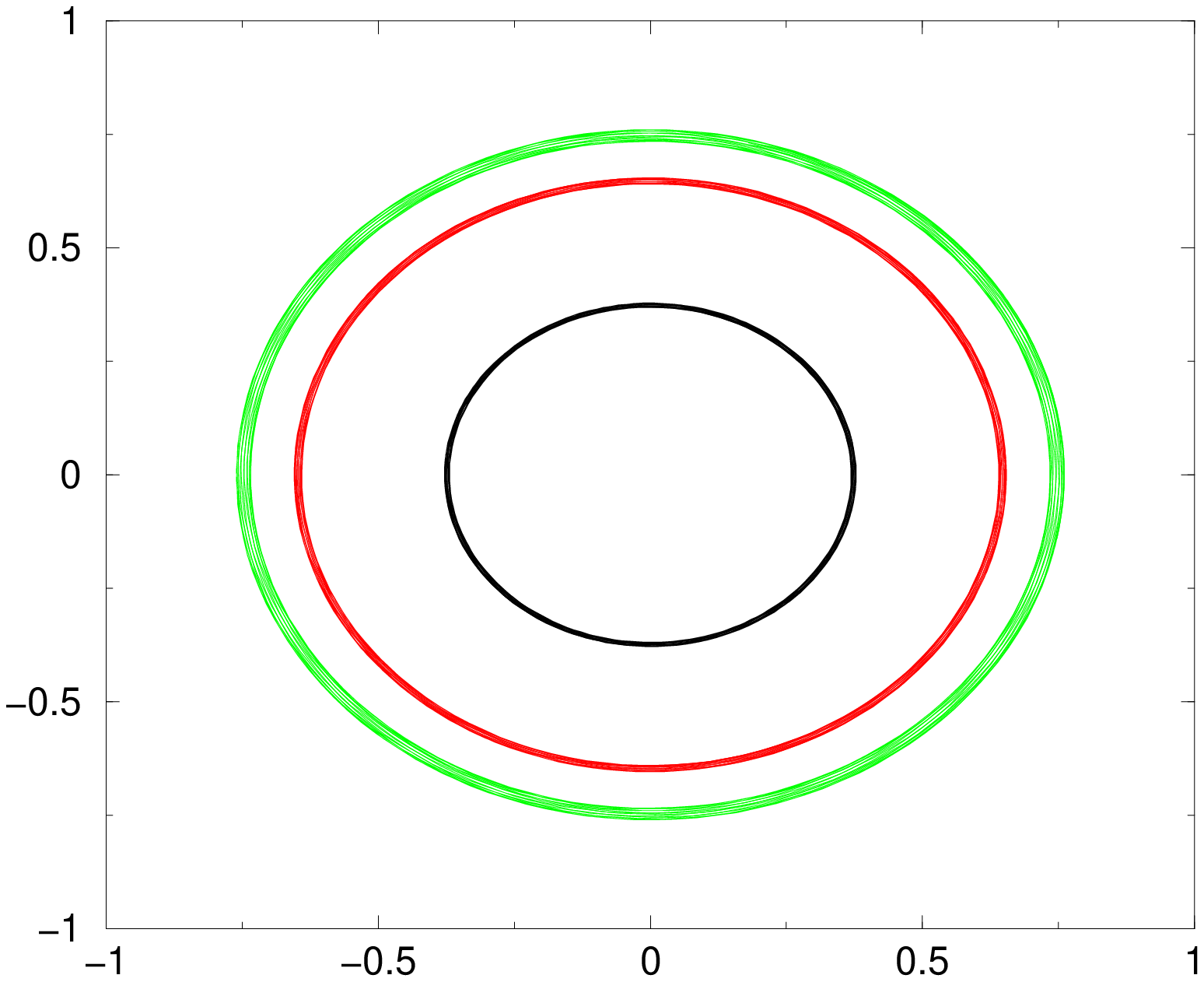} \hspace{5mm}
\includegraphics[angle=0,height=60mm,draft=false]
{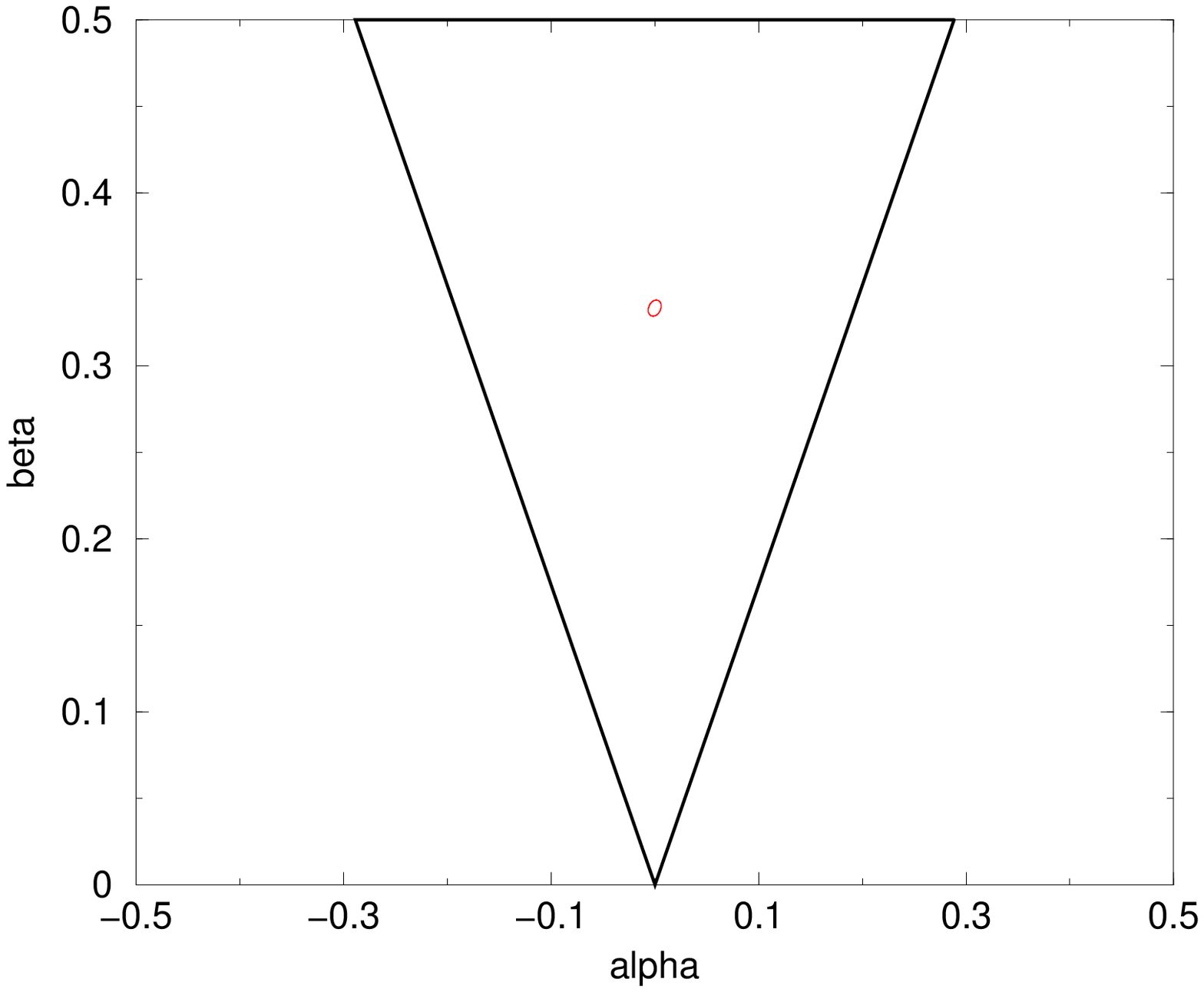}
}
\caption{Trajectory of the three vortices with respect to the centroid (left)
and in the $\alpha-\beta$ plane (right): case 1 (top), case 4 (middle)
and case 7 (bottom).  The trajectory of vortex 1 is depicted in black, while
those of vortices 2 and 3 are shown in red and green, respectively.}
\label{fig:traj}
\end{figure}

\begin{figure}[phtb]
\centerline{
\includegraphics[angle=0,height=42mm,draft=false]
{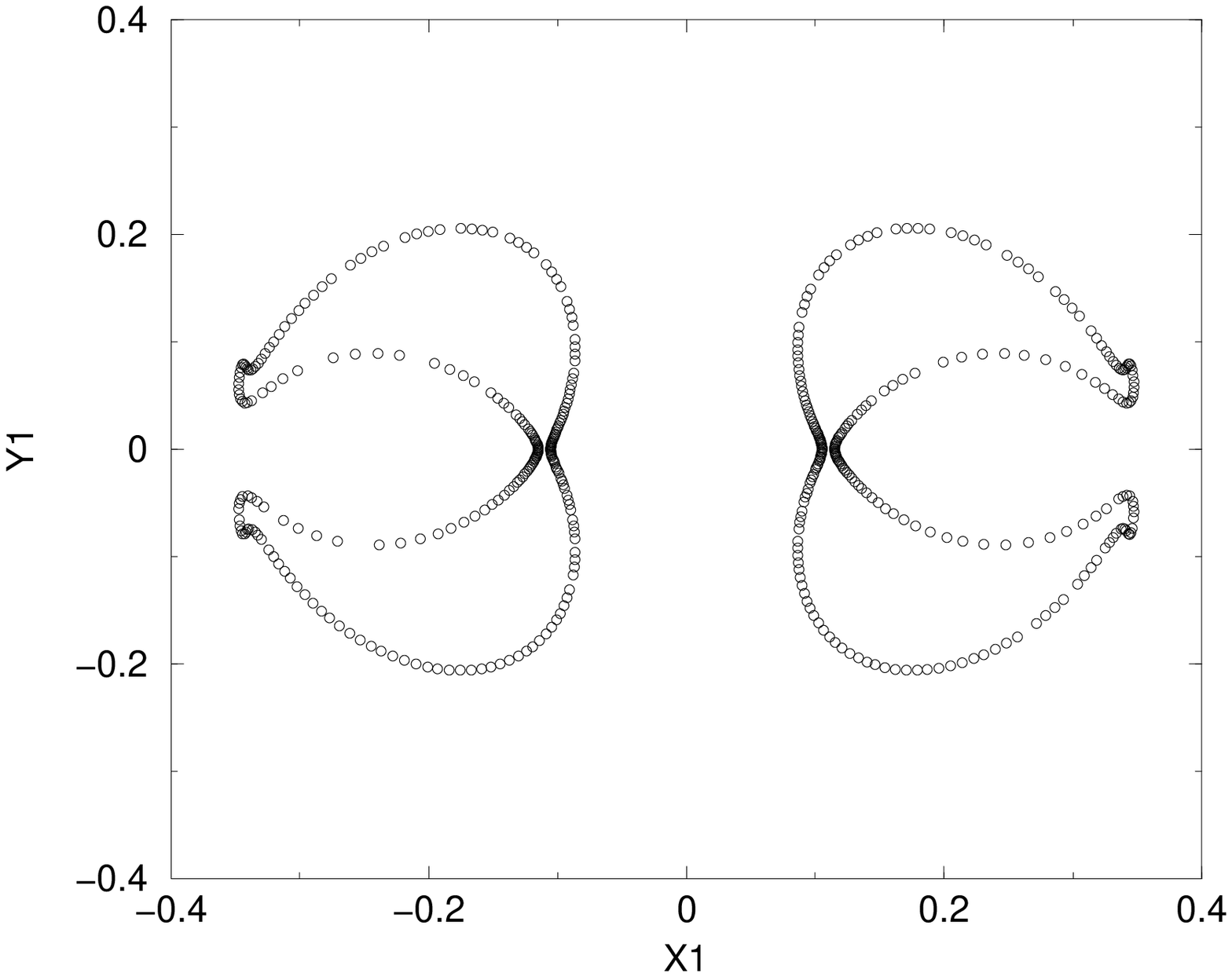} \hspace{5mm}
\includegraphics[angle=0,height=42mm,draft=false]
{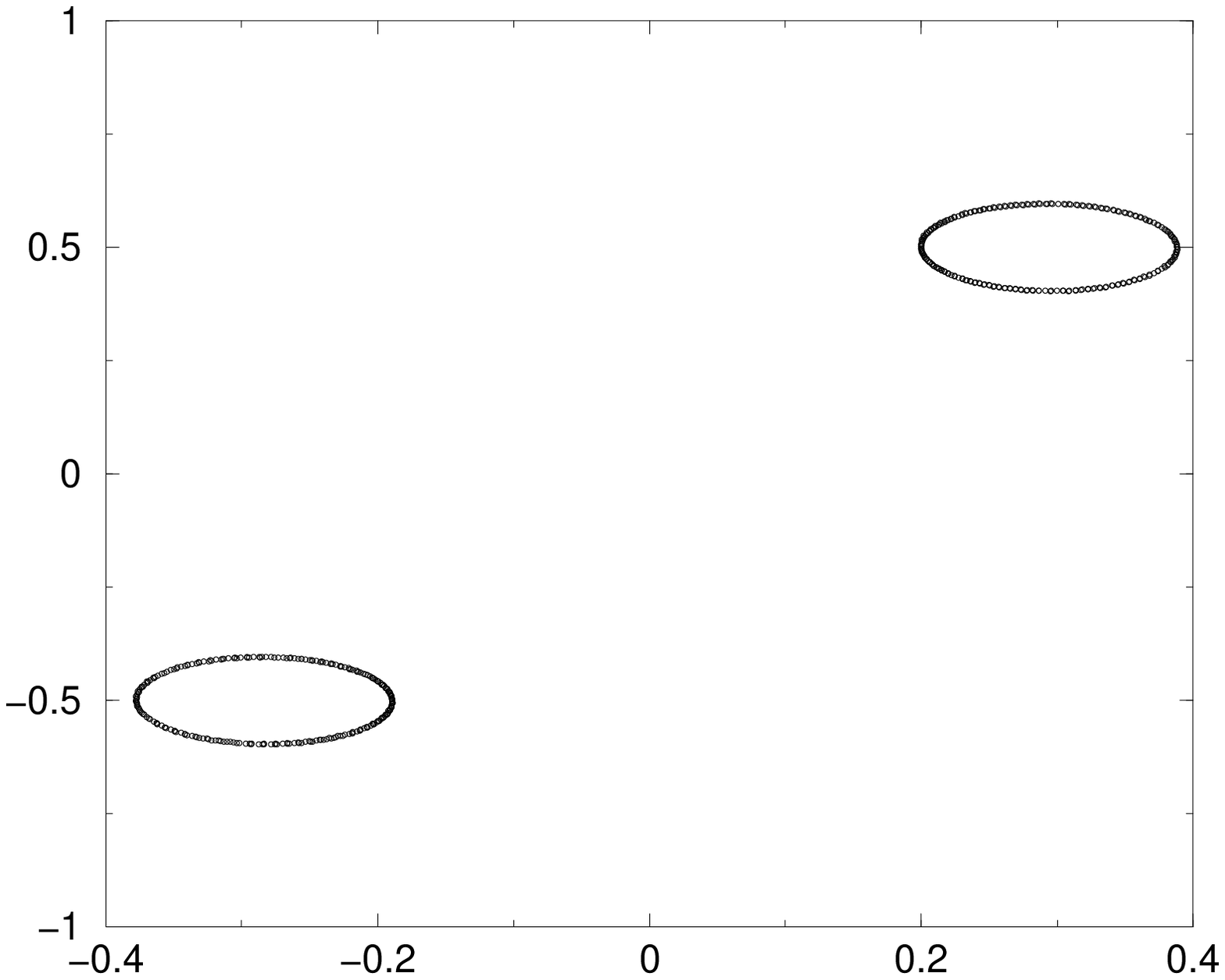} \hspace{5mm}
\includegraphics[angle=0,height=42mm,draft=false]
{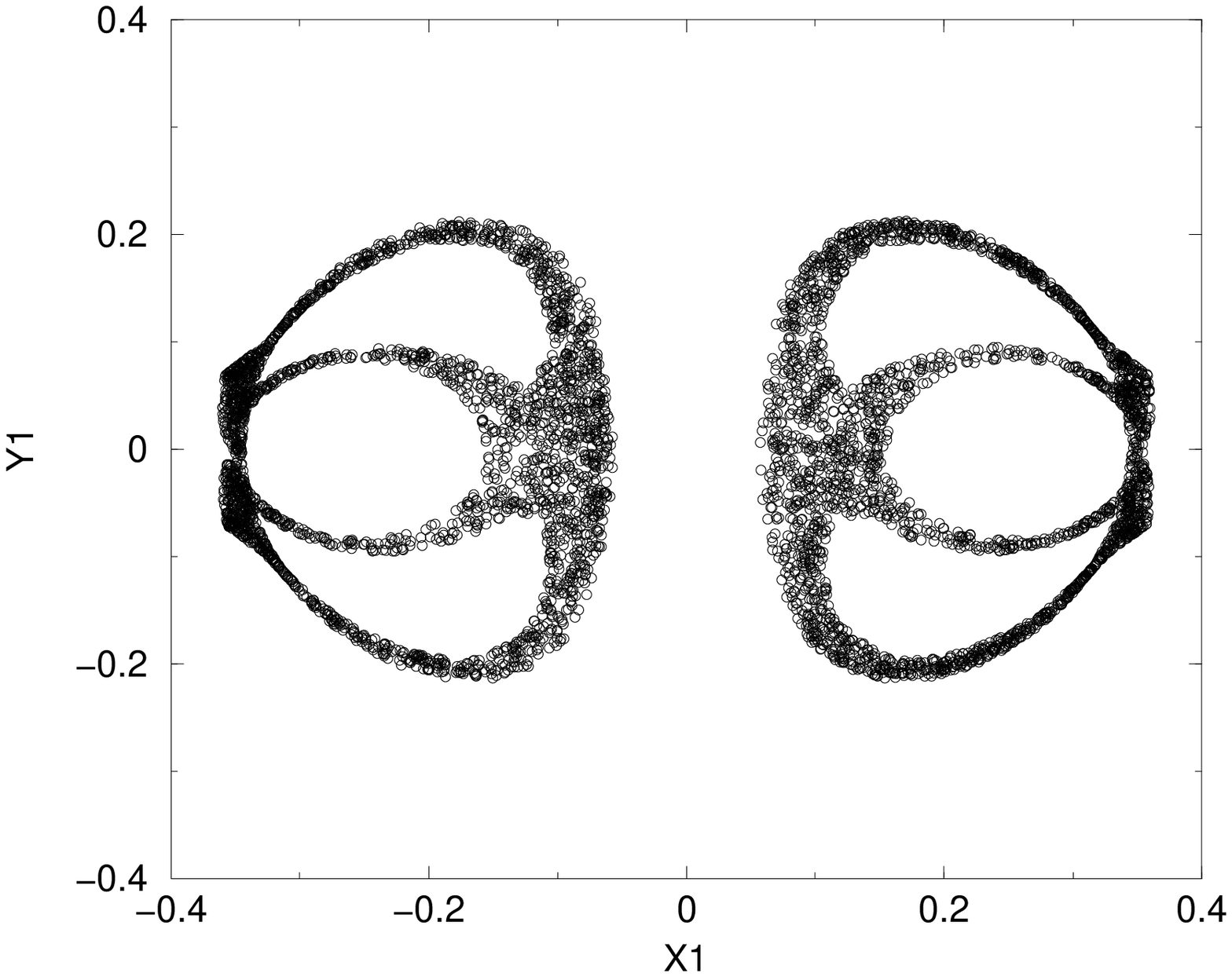} } \caption{Poincar\'e maps for case 4.  Left: unperturbed
system; middle: type A perturbation with $h=10$; right: type A
perturbation with $h=2.5$.} \label{fig:case4}
\end{figure}

\begin{figure}[phtb]
\centerline{
\includegraphics[angle=0,height=42mm,draft=false]
{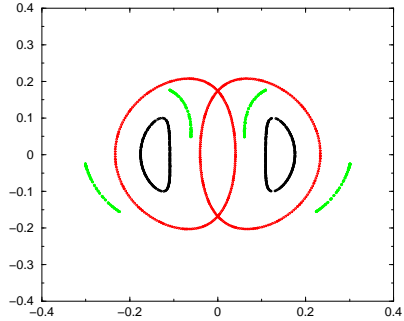} } \caption{Poincar\'e maps for case 1 (black),  4
(red), and 7 (green). Type C perturbation with $\rho = 250$.}
\label{fig:typII}
\end{figure}

\begin{figure}[phtb]
\centerline{
\includegraphics[angle=0,height=42mm,draft=false]
{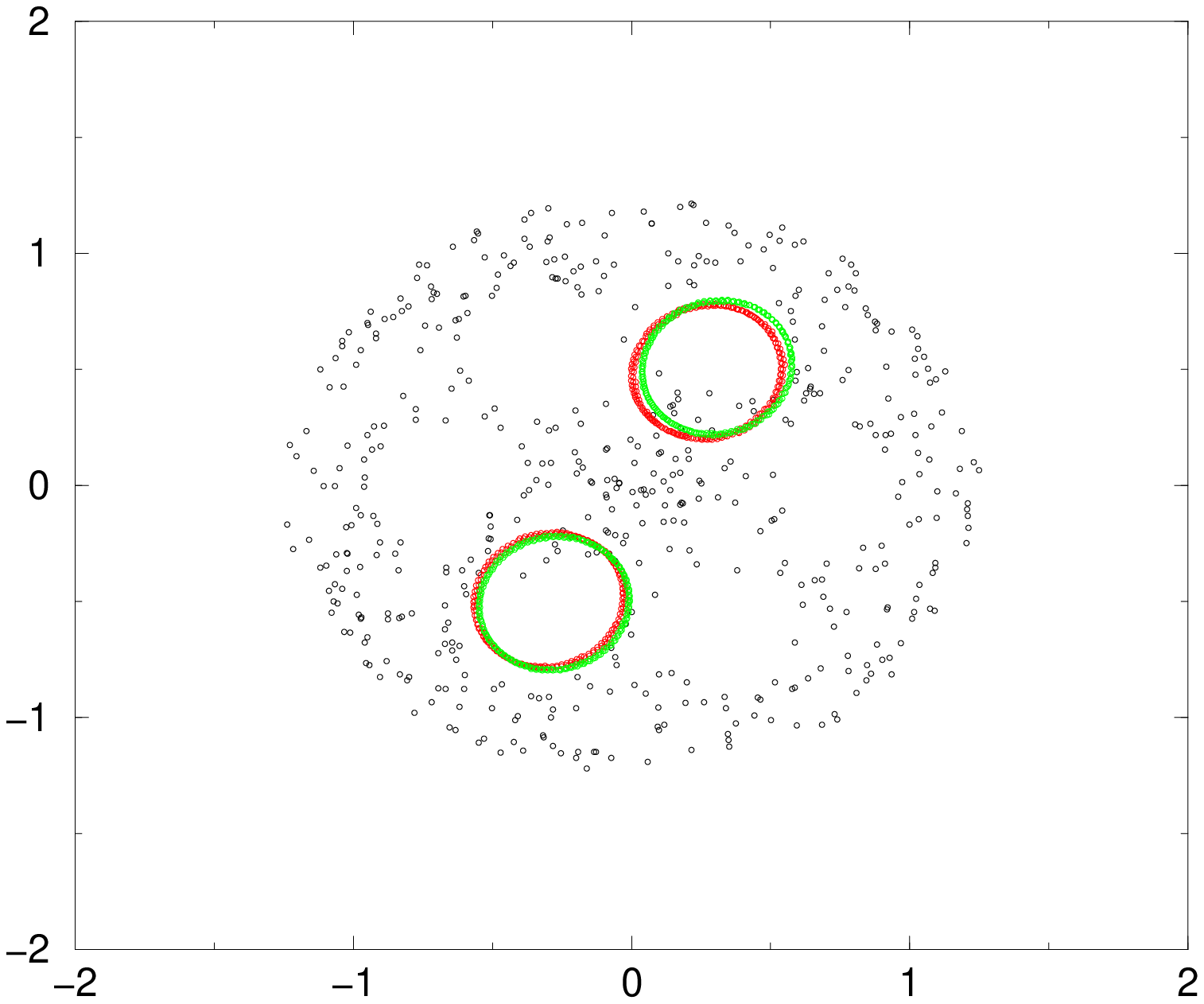} \hspace{5mm}
\includegraphics[angle=0,height=42mm,draft=false]
{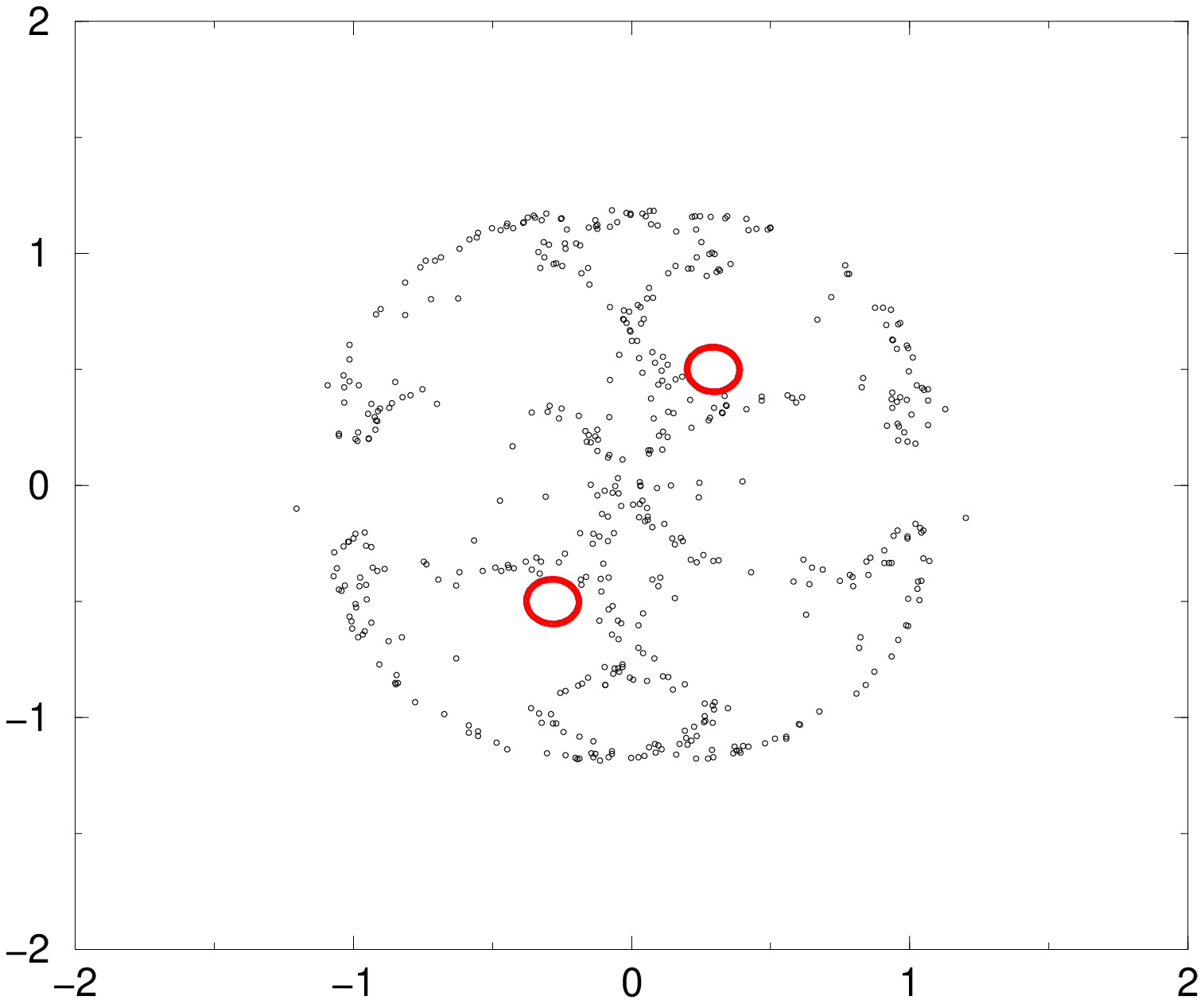} } \caption{Poincar\'e maps for the restricted 4-vortex
configuration.  Left: maps for case i (black), iii (red) and v
(green). Right: maps for case ii (black) and iv (green). }
\label{fig:4vortex}
\end{figure}

\begin{figure}[phtb]
\centerline{
\includegraphics[angle=0,height=60mm,draft=false]
{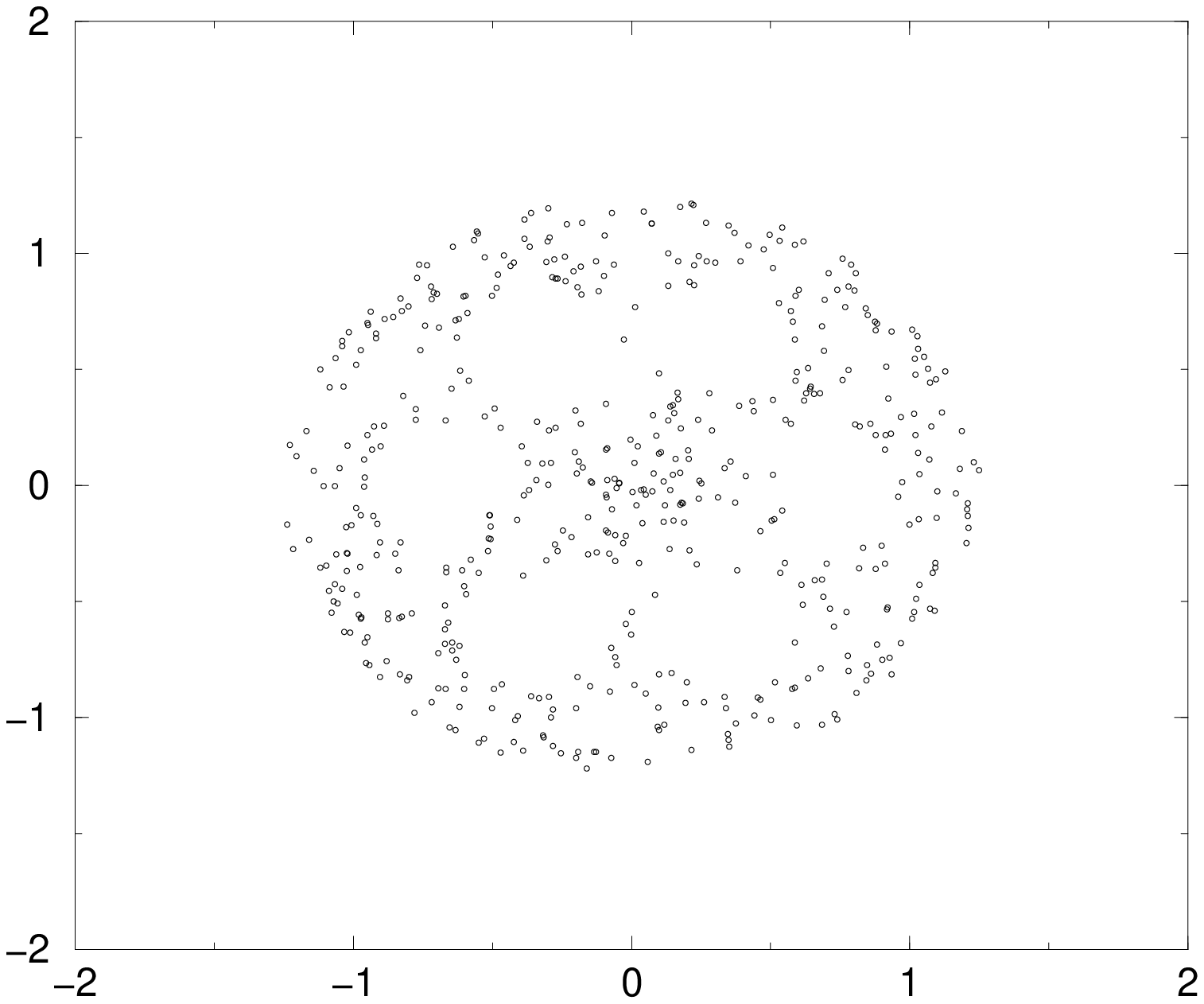} \hspace{5mm}
\includegraphics[angle=0,height=60mm,draft=false]
{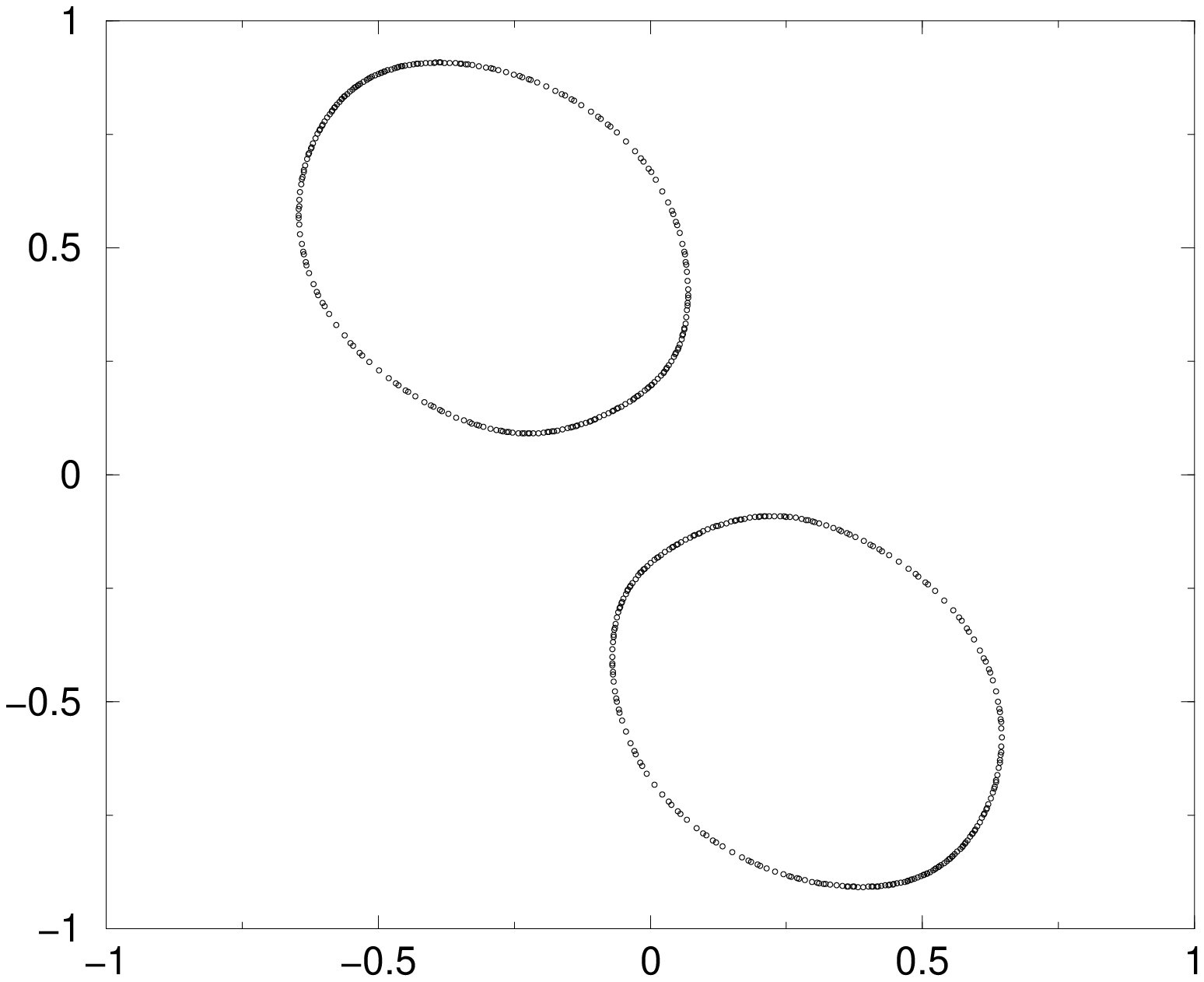} } \centerline{
\includegraphics[angle=0,height=60mm,draft=false]
{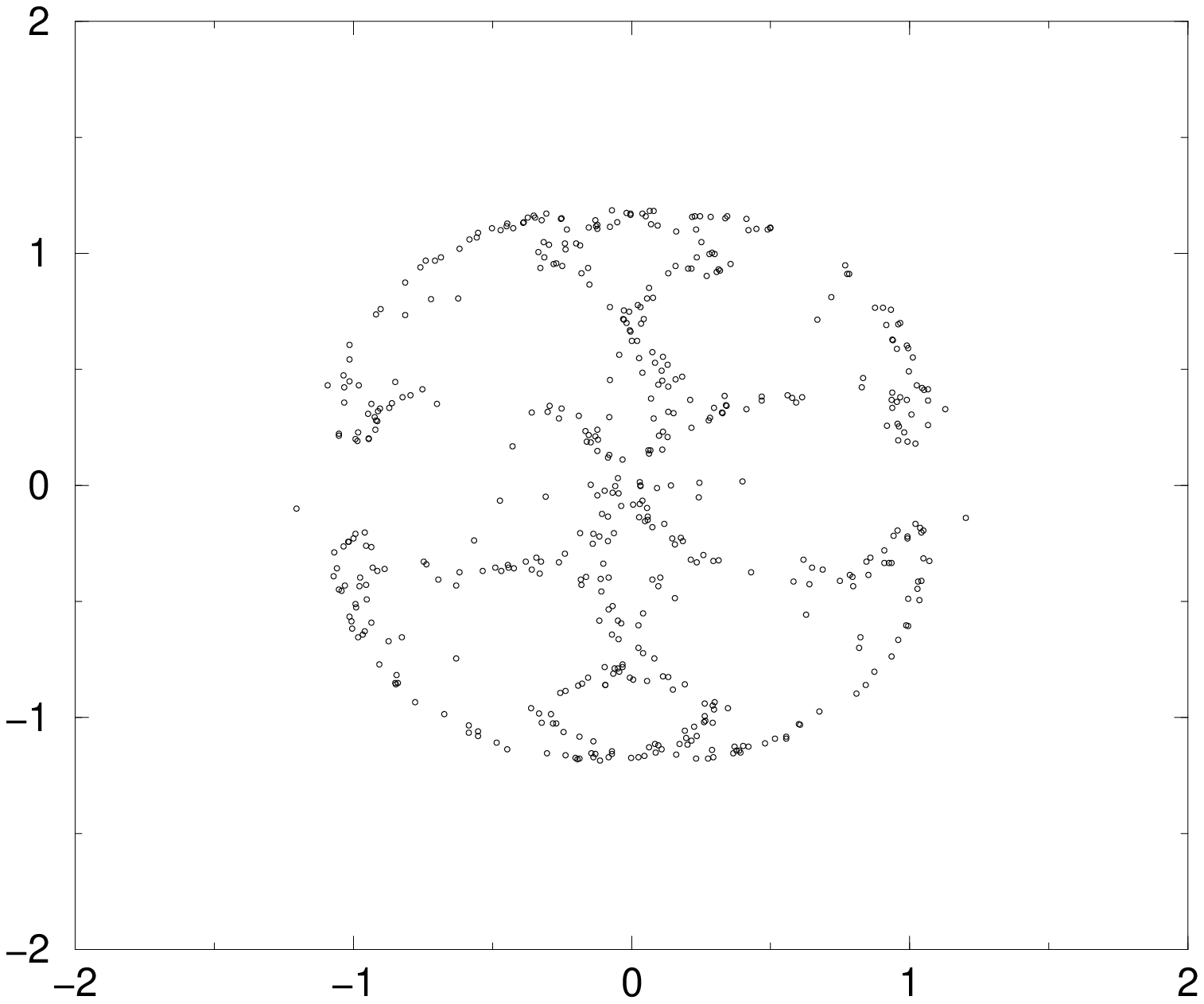} \hspace{5mm}
\includegraphics[angle=0,height=60mm,draft=false]
{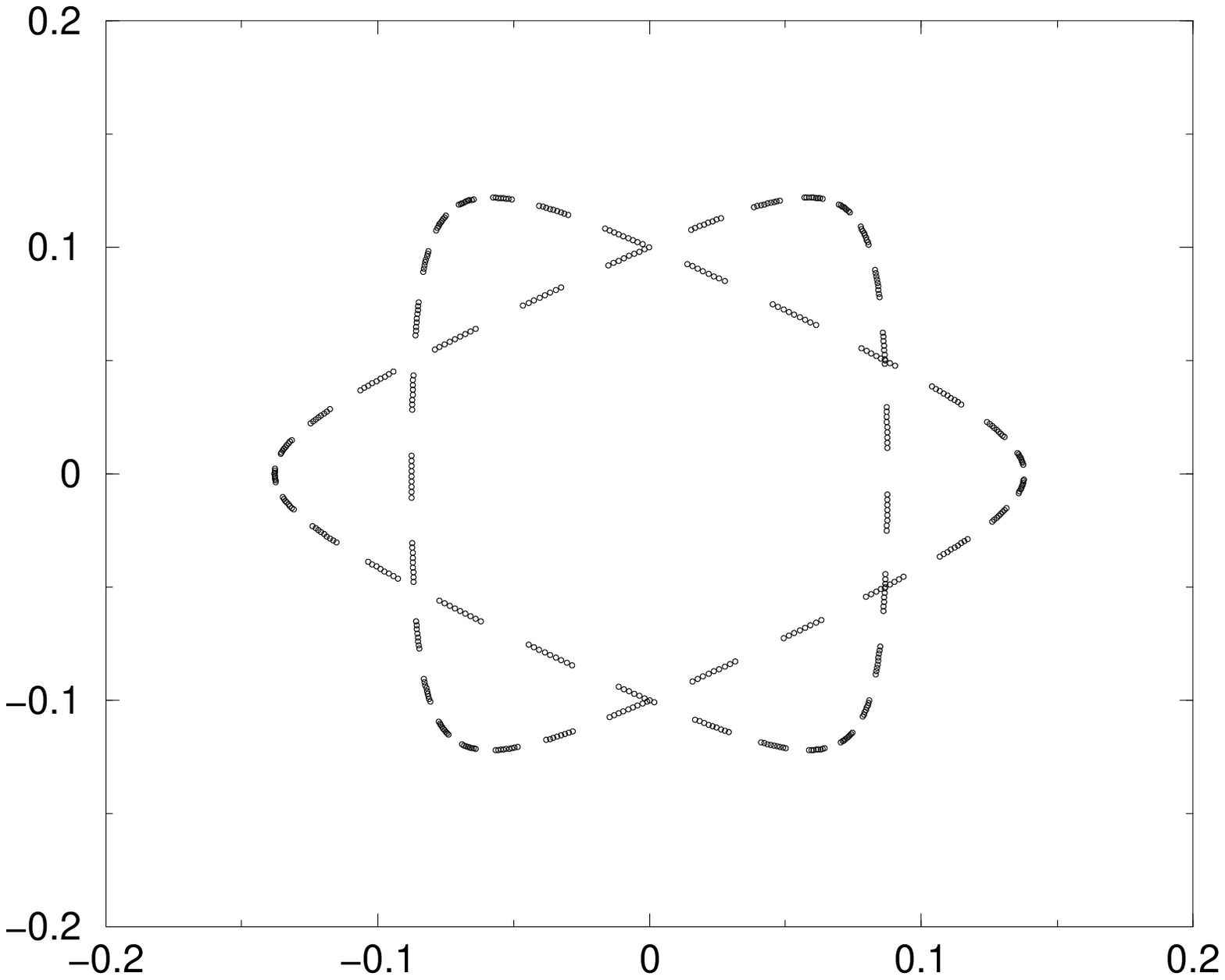} } \centerline{
\includegraphics[angle=0,height=60mm,draft=false]
{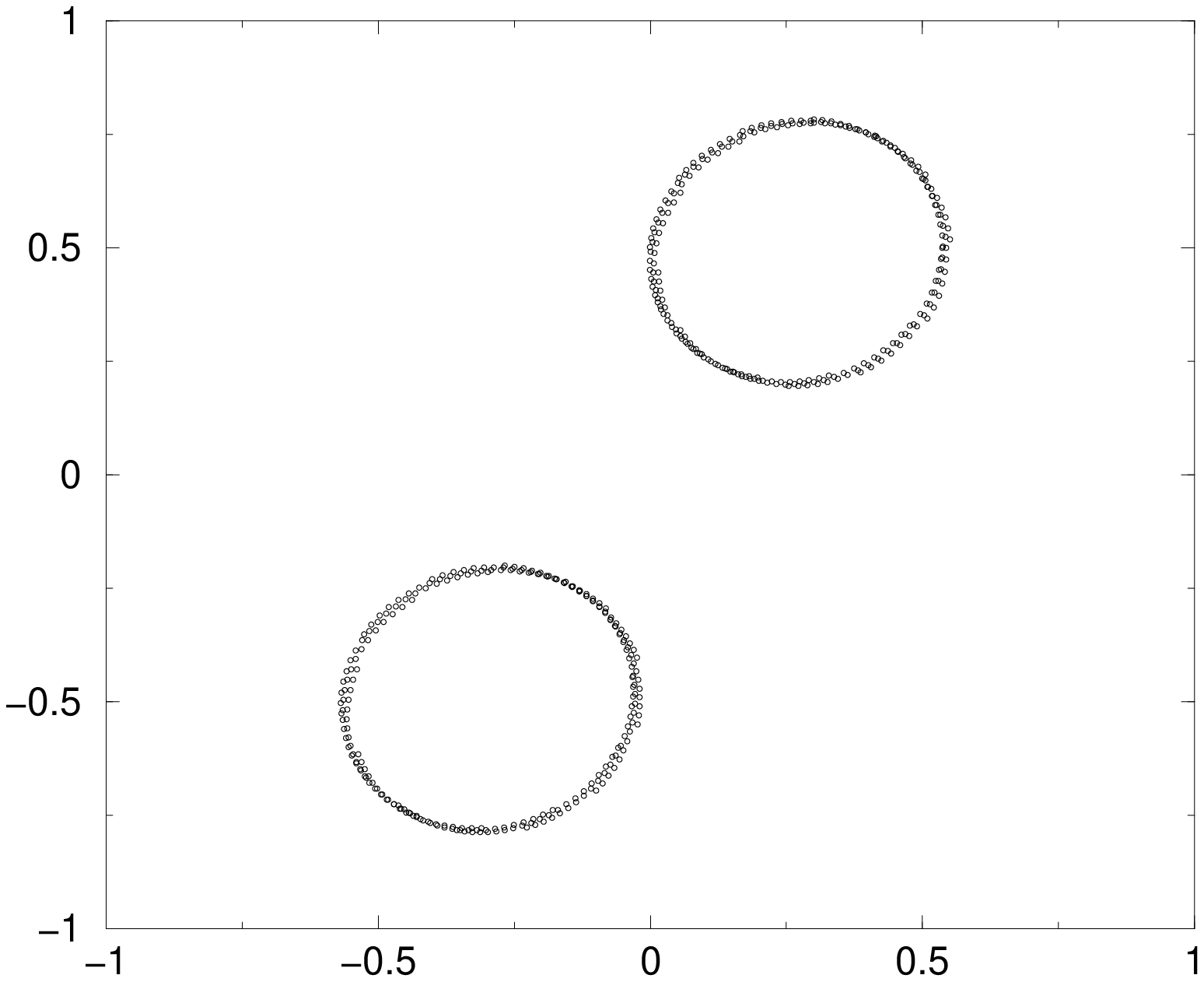} \hspace{5mm}
\includegraphics[angle=0,height=60mm,draft=false]
{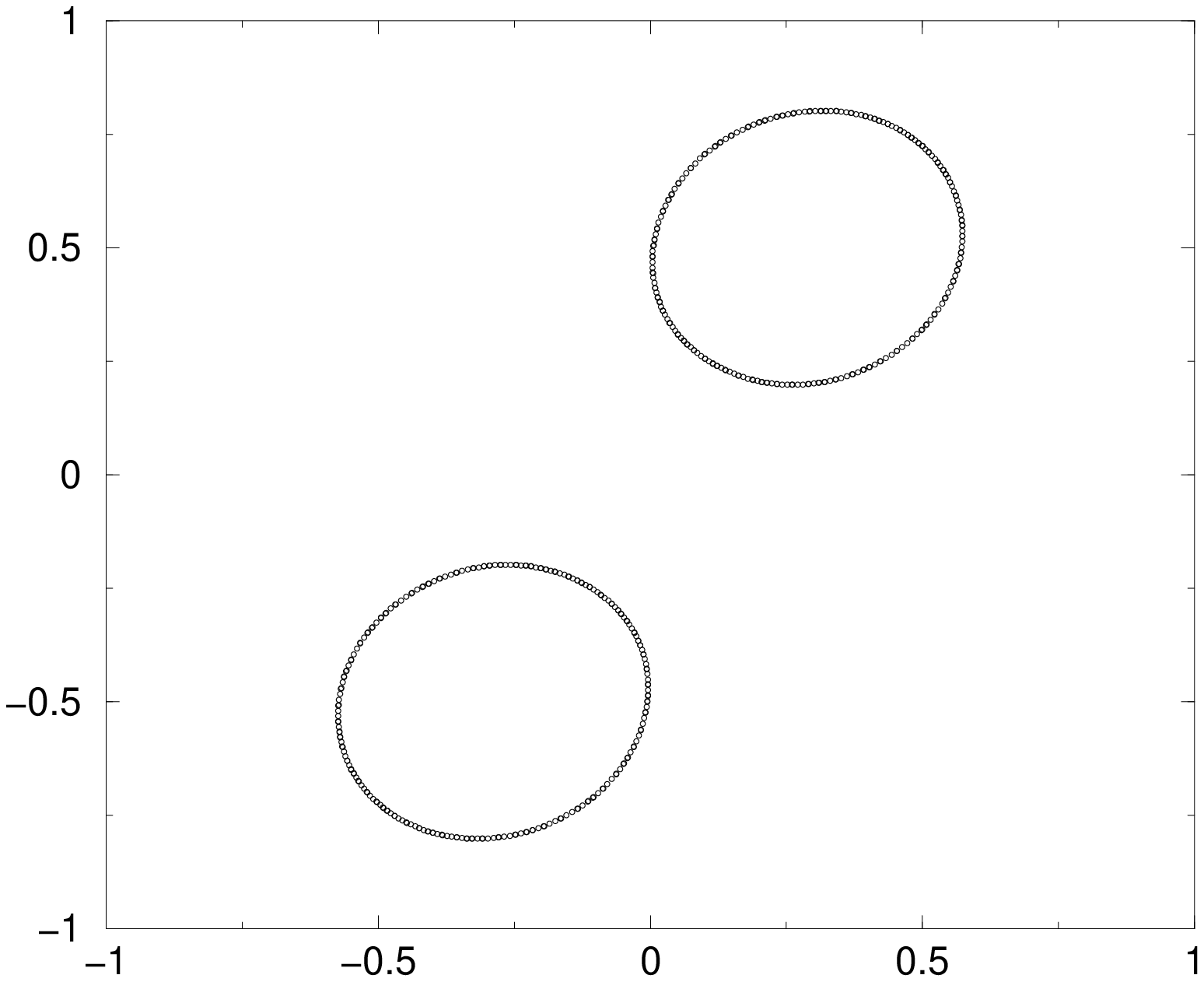} } \centerline{
\includegraphics[angle=0,height=60mm,draft=false]
{p4.eps} \hspace{5mm}
\includegraphics[angle=0,height=60mm,draft=false]
{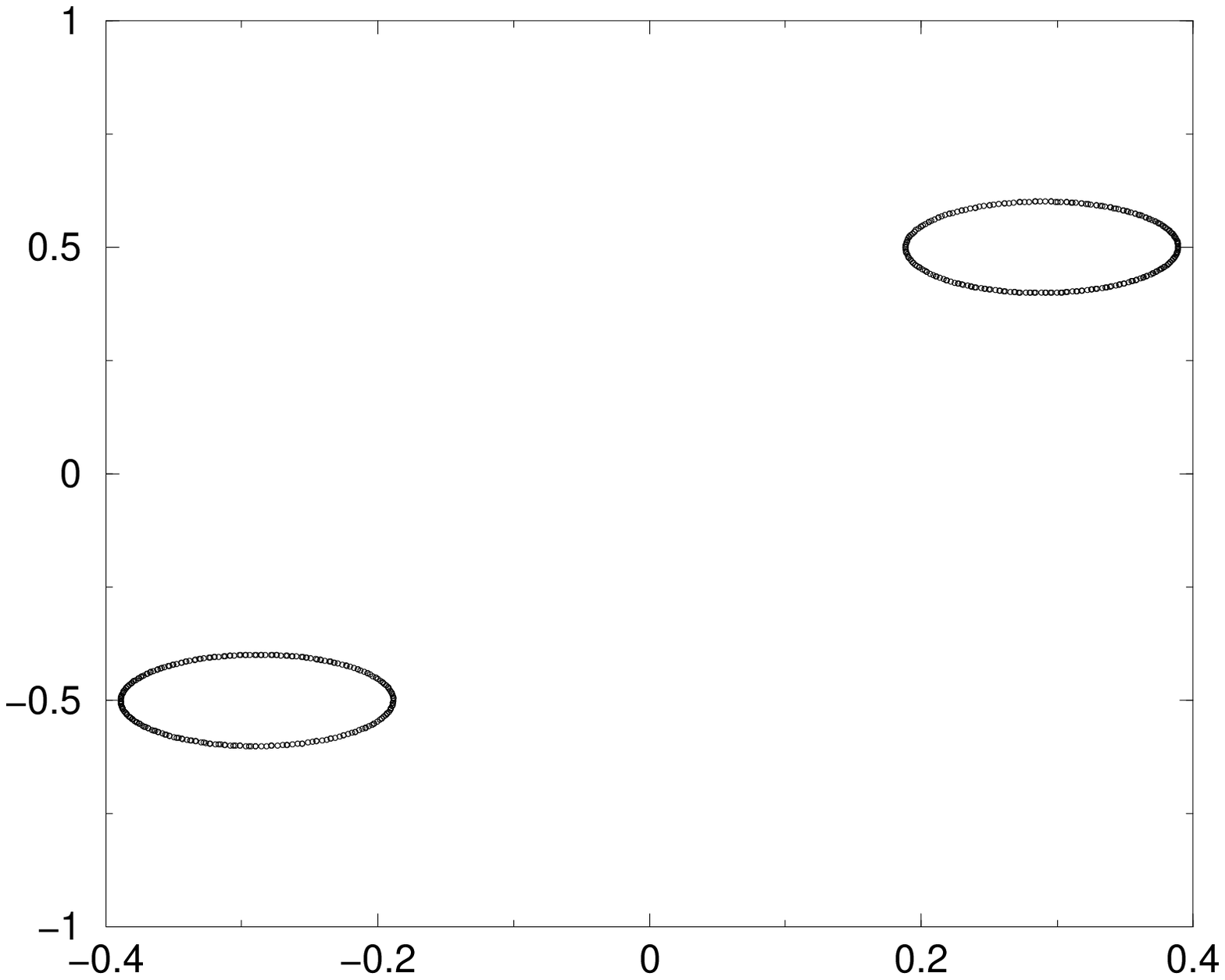} } \caption{Poincar\'e maps for the restricted 4-vortex
configuration, for perturbed (left) and unperturbed (right) systems.
Shown are results for cases i, ii, iii, and iv, arranged from top.}
\label{fig:4vmaps}
\end{figure}

\noindent\textsf{\textbf{Summary of Simulation Results}}
\begin{itemize}
\item Our main focus is on examples where the unperturbed system is of the elliptic type, \emph{i.e.}
cases 1, 4, 7.  Initial configuration near $Q_3$, in between two
separatrices (region 3B shown in Fig. 1), and near the center $E$.
We expect more sensitivity to perturbations for of any of the three
types for case 4 than 1 and 7, since case 4 represents a
configuration closer to the separatrix shown in Fig. 1, and this is
borne out by our simulations.
\item The motion of the three vortices in the plane for cases 1, 4, and 7 in the unperturbed system is
shown in Figure~\ref{fig:traj}. Notice that even in this case, the
trajectories are considerably more complicated for case 4 than cases
1 and 7. Again this is due to proximity to the separatrix.
\item For a type A perturbation with $h=10$, we observe little impact on trajectories of the 3 vortices.
This is consistent our theoretical results.
\item For cases 1 and 7, we also consider lower values of $h$ - as low as $h=2.5$.  For case 4 corresponding,
to a point close to the separatrix for the unperturbed system, the
dynamics become much more complicated at a more rapid pace as $h$
decreases. The Poincar\'e maps ($x_1,y_1$ for $y_2 = 0$) shown in
Fig.~\ref{fig:case4} suggest that when $h=2.5$ the dynamics is
starting to exhibit chaotic regimes.
\item Also considered in our simulations are type C perturbations. We found that it is necessary to start
with a large mean radius for the three slender coaxial vortex rings
in order to observe theoretical predictions. With $\rho=250$, the
motion remains regular. This is illustrated in the Poincar\'e maps
shown in Fig.~\ref{fig:typII}.
\item In the final set of simulations, we considered examples of restricted four vortex dynamics
The configuration of the three larger (principal) vortices, which we
designate as vortices 1-3, is defined in terms of equal-strength
vortices, $k_i = 1$ for $i=1,2,3$, initially located on the vertices
of an equilateral triangle of side 1. Specifically, the initial
conditions are given by: $z1=(-0.5,0)$, $z2=(0.5,0)$, and $z_3 =
(0,h)$ where $h=\sqrt{3}/2$.  For the unperturbed case, the
trajectories with respect to the centroid are identical circles
centered at the origin. We consider perturbation due to the
introduction of a weak vortex having $k_4 = 0.1$.  Five different
initial conditions for the fourth vortex are considered, namely $z_4
= (0.1,0)$, $z_4 = (0.1,h/3)$, $z_4 = (0.1,2h/3)$, $z_4 = (0.1,h)$,
and $z_4 = (0.1,4h/3)$.  We refer to these as cases i - v,
respectively.
\item Consistent with predictions, the trajectories of vortices 1-3 are only weakly affected by the
introduction of the fourth (weaker or smaller) vortex (not shown).
More precisely, we find that the common circle of motion of the
three principal vortices becomes slightly thickened (forming a thin
annulus) in the plane - with the thickness of the annulus and
variation of motion within it reflecting the magnitude (as well as
the particular form) of the perturbation caused by the smaller
vortex. One expects the annulus to thicken when the smaller vortex
starts closer to the configuration of vortices 1-3, its strength
increases, or when the weaker vortex starts in a position that has a
more subtle connection with the configuration of the whole set of
four vortices. These subtleties will be discussed at some length in
the following points.
\item The most dramatic effect on the dynamics for the restricted four vortex problem is naturally
manifested in the motion of the fourth (smaller or weaker) vortex.
The trajectory of the fourth vortex can be either regular or complex
and even chaotic, depending on its initial location and strength.
Moreover, the complexity of its trajectory is quite sensitive to
changes in its initial location and strength. This is illustrated in
the Poincar\'e maps of  ($x_4,y_4$ at $y_1=0$).
\item The Poincar\'{e} maps for the motion of the fourth vortex shown in Fig.~\ref{fig:4vortex}
indicate via the characteristic splattering that the dynamics for
cases i and ii exhibit typical chaotic behavior, or at least a
transition from regular to chaotic motion, whereas cases iii and iv
are essentially regular. This may appear to be somewhat surprising
in light of the fact that the distance from vortex 4 to the set of
vortices 1-3 in cases i and ii is not very different from the
distances in cases iii and iv.
\item A plausible explanation for the somewhat counter intuitive
behavior described above is as follows: If we view vortex 4 as
perturbing the motion of vortices 1-3 (see Fig. 1), we would not, as
indicated above, expect this to have much effect on the dynamics of
vortices 1-3. In particular, the configuration of vortices 1-3
places it initially at the center $E$, so the inherent stability
should preserve periodic motion for this configuration that varies
only slightly from an equilateral array. Dually, we can consider the
triple of vortices 1, 2 and 4, as being perturbed by the strong
third vortex. For cases i and ii, the initial configuration of
$\{1,2,4\}$ is rather close to a separatrix shown in
Fig.~\ref{fig:4vortex}. This experiences a strong perturbation from
vortex 3, which is substantial enough to push the configuration
across the separatrix and also possibly break this curve. Similarly,
cases iii-v are such that the initial configuration of vortices i,
ii and iv are substantially further removed from such a separatrix,
so then it is not surprising that the dynamics for cases i and ii is
far less regular than that for cases iii and iv.
\item Some of these subtle points are further illustrated and
contrasted in Fig.~\ref{fig:4vmaps}, which shows the dynamics of the
fourth vortex, again using Poincar\'e maps. The dynamics of the
fourth vortex of finite strength, portrayed in the left column of
figures, is contrasted with the same cases on the right, where the
strength of vortex 4 is taken as zero. In other words, we are
treating vortex 4 as a passive particle propelled by the motion of
the three stronger vortices. Notice that in all cases where the
strength of the vortex is zero, the Poincar\'e maps show that the
motion of this vortex is regular. The juxtaposition in
Fig.~\ref{fig:4vmaps} underscores the point that the effect of the
fourth vortex is dependent both on its initial placement and
strength.

\end{itemize}

\clearpage

\bigskip

\noindent\textsf{\textbf{VII. CONCLUDING\ REMARKS}}

\medskip

Following a brief summary of the dynamics of three point vortices
moving in an ideal fluid in the plane in both a Hamiltonian and
trilinear coordinate context, which highlighted the properties most
pertinent for our investigation - such as the existence of a center
regardless whether the unperturbed system is of elliptic, hyperbolic
or parabolic type corresponding, respectively, to $\kappa_1\kappa_2
+ \kappa_1\kappa_3 + \kappa_2\kappa_3$ positive, negative or zero
(assuming that $\kappa_1+\kappa_2+\kappa_3 \neq 0$) - we described
the three kinds of perturbations that comprised the focus of this
paper. These were three vortices in a half-plane, a restricted four
vortex problem, and the motion of three coaxial vortex rings.

We then formulated and proved some new results on the existence of regular
regimes - characterized by behavior associated to integrable Hamiltonian
systems - for non-integrable perturbations of LA-integrable three vortex
dynamics. The regular dynamics in point includes the existence of an ample set
of quasiperiodic flows on invariant KAM tori along with periodic solutions.
These results, which were proved employing rather modern variants and
generalizations of the KAM and Poincar\'{e}-Birkhoff fixed point theorems,
extend some related recent results for point and coaxial ring vortices to
cases where the signs of the vortex strengths differ, and are such that their
sufficient conditions can be readily checked. But, they are really just
existence theorems that provide little information about the location and
finer geometric details of the predicted invariant tori and periodic orbits.

The lack of quantitative specificity in our first set of perturbation theorems
was subsequently addressed for the case of periodic orbits, with some novel
results concerning the existence of periodic solutions for perturbed systems
that closely approximate cycles of the unperturbed three vortex system. These
results - possessing the look and feel of classical perturbation theorems -
were proved using Poincar\'{e} sections in concert with the Brouwer fixed
point theorem, are stated in terms of easily verifiable criteria guaranteeing
the existence of periodic orbits of the perturbed system close to those of the
three vortex system.

Several numerical simulations were performed to illustrate the perturbations
theorems obtained in a variety of contexts. Also included were examples to
indicate how the desired regularity of the dynamics breaks down as the
assumptions of our theorems are stretched to the limit, which is accomplished
by showing how there is a (transitional) preponderance of chaotic regimes as
these limits are exceeded.

Our study of perturbations of the three vortex problem has, we hope,
shed some new light on these particular systems as well as more
general vortex phenomena. It certainly has - as is often the case in
such studies - created more interesting questions than it has
answered, several of which we intend to investigate in the near
future. For example, our work here suggests that most of our results
can be extended to systems comprised of a larger number of vortex
elements if they are clustered in certain ways. As an illustration
of this, suppose that one has a pair of vortices of strengths of the
same sign that are initially quite close together. Then this pair
constitutes a binary system that should, under the right conditions,
remain close together for all time, and behave as if it were a
single vortex of strength equal to the sum of the strengths of each
of its components if any other point vortices have initial distances
from this pair that are substantially larger than the distance
between the pair. A configuration of four vortices including such a
binary system and two other vortices initially quite distant from
the pair, might well be expected to behave very much like a slightly
perturbed three vortex system in the plane (with the binary system
acting like a single vortex), thus generating a profusion of regimes
exhibiting regular (integrable type) dynamics. Interesting
clustering problems of this kind are something that we plan to turn
our attention to in our future research.

Another rather manifest question raised by our investigation - that we intend
to presently turn our attention to - concerns extending the quantitative
perturbation results obtained for periodic orbits to invariant tori. It
appears that it may be possible, using some type of generalization of the
Poincar\'{e}-Birkhoff fixed point theorem, to show that in certain cases the
perturbed system has an invariant torus that is close to an invariant torus
for the unperturbed system.

\bigskip

\noindent\textbf{\textsf{ACKNOWLEDGMENTS}}

\medskip

The authors wish to dedicate this paper to Professor Theodore Yaotsu
Wu of the California Institute of Technology, belatedly on his
eightieth birthday, and to honor his life-long pioneering
contributions to fluid mechanics and soliton waves, and his
leadership in training younger scientists and engineers. In
particular, L. T. would like to refresh his memories with Professor
Wu of being classmates in the ninth grade at Guang Hua Middle School
and 1946 alumni of the Chioa Tung University, Shanghai.

In addition, the authors wish to express their appreciation to
Charles Doering and Paul Newton, guest editors of the special focus
issue of the Journal of Mathematical Physics on mathematical fluid
dynamics, for their gracious invitation to contribute this paper.

Lastly, we wish to acknowledge that preliminary versions of some of
the material in this paper were reported on at the ICTAM Conference
in Warsaw, Poland, 2004, the Third MIT Conference on Computational
Fluid and Solid Mechanics, 2005, and the GAMM Annual Meeting,
Berlin, 2006, and appeared in summary form in the proceedings of the
first two conferences.

\pagebreak

\end{document}